\documentclass[12pt]{article}

\usepackage{tipa}
\usepackage{amsfonts,amssymb}
\usepackage{amssymb}
\usepackage{float}
\usepackage{graphicx}
\usepackage{amsthm}
\usepackage{amsmath}
\usepackage[hang]{caption2}%\usepackage{caption2}
\usepackage[all]{xy}\xyoption{all}

\newtheorem{lemma}{Lemma}
\newtheorem{theorem}{Theorem}

\newtheorem{remark}{Remark}
\newtheorem{definition}{Definition}

\begin{document}

\title{Some Symmetric Orbits for N-body Type Difference Equations}

\author{Leshun Xu$^{1,2}$, Yong Li$^{1}$ and Menglong Su$^{1}$\\
              $^{1}$College of Mathematics,\\
              Jilin University, Changchun 130012, P.R. China \\
            $^{2}$\emph{Present address:} College of Mathematics and Computer Science,\\
            Nanjing Normal University, Nanjing
            210097, P.R. China.     \\
e-mail: lsxu@sohu.com; liyong@jlu.edu.cn; mlsu@email.jlu.edu.cn
}
\date{25.10.2006}
\maketitle

\begin{abstract}
This paper introduces a new difference scheme to the difference
equations for N-body type problems. To find the non-collision
periodic solutions and generalized periodic solutions in
multi-radial symmetric constraint for the N-body type difference
equations, the variational approach and the method of minimizing
the Lagrangian action are adopted and the strong force condition
is considered correspondingly, which is an efficient method in
studying those with singular potentials. And the difference
equation can also be taken into consideration of other periodic
solutions with symmetric or choreographic constraint in further
studies.

\end{abstract}

\section{Introduction}\label{intro}
The N-body problem is a classical and important problem in celestial
mechanics and mathematics. It consists of determining the orbits of
N bodies interacting in accordance with the gravitational law of
Newton. That is two bodies attracting each other, the force of
attraction being directed along the line joining them, proportional
to the product of the masses, and inversely proportional to the
square of the distance between them. For given $N$ bodies, the
preceding N-body problem can be described by the following nonlinear
system of second order differential equations:
\begin{equation}\label{Nb}
-m_{i}\ddot{q_{i}}=\sum_{l=1,l\neq
i}^{N}\frac{Gm_{i}m_{l}(q_{i}-q_{l})}{|q_{i}-q_{l}|^{3}},\quad\quad
i=1,\cdots,N,
\end{equation}
where $m_{i}$ is the mass of $i$th body, $q_{i}$ is the position of
$i$th body in $R^{d}$, $N$ is the number of bodies, and $G$ is the
universal gravitational constant, which is always taken to be 1 for
convenience in mathematics. Hence we take $G=1$ in this paper.

\begin{figure}[H]
    \begin{minipage}[t]{0.3\linewidth}
%    \centering
    \includegraphics[width=1.6in]{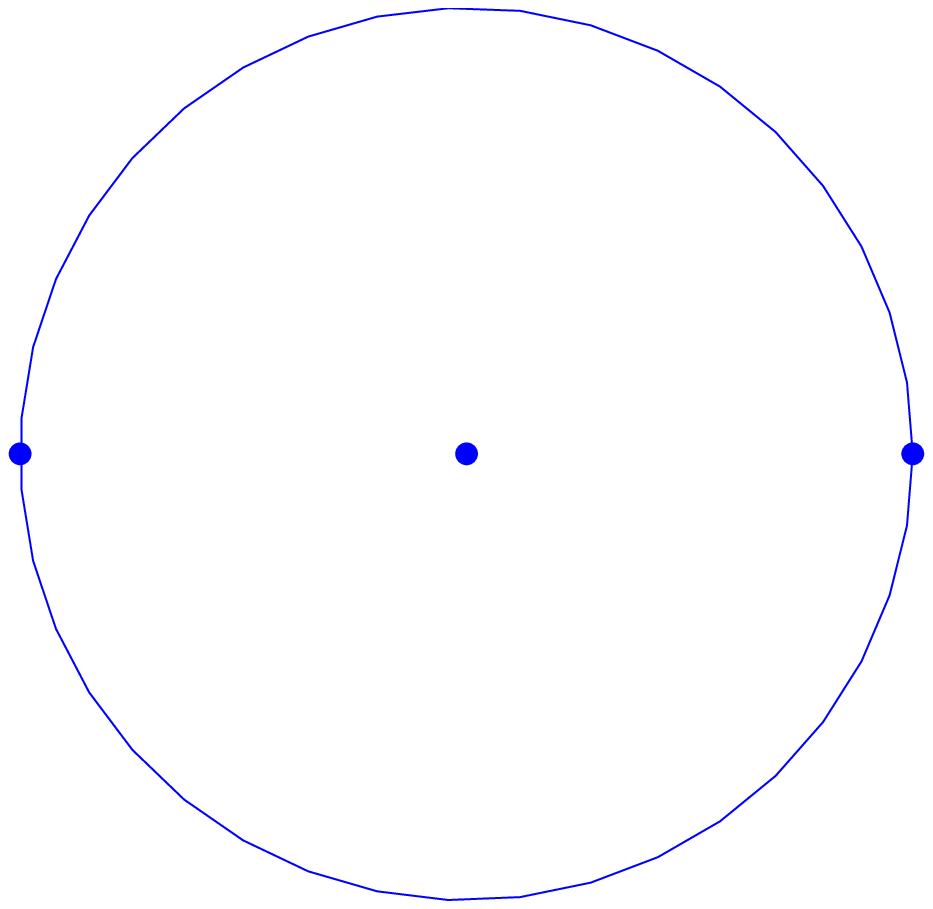}
%    \parbox{1.6in}{\caption{$m=[1,1,1]$,}\label{N_104_0126}}
    \caption{Euler's Solution}\label{Euler}
    \end{minipage}
    \hspace{0.5cm}
    \begin{minipage}[t]{0.3\linewidth}
%    \centering
    \includegraphics[width=1.6in]{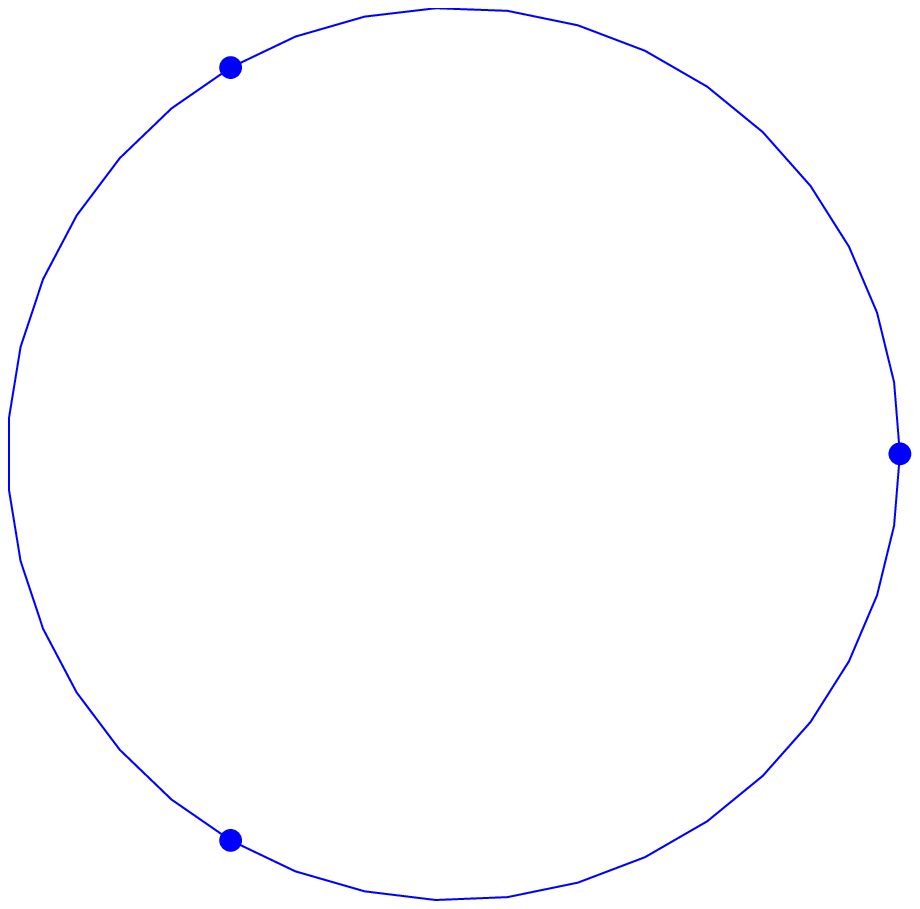}
    \caption{Lagrange's Solution}\label{Lagrange}
    \end{minipage}
    \hspace{0.5cm}
    \begin{minipage}[t]{0.3\linewidth}
%    \centering
    \includegraphics[width=1.6in]{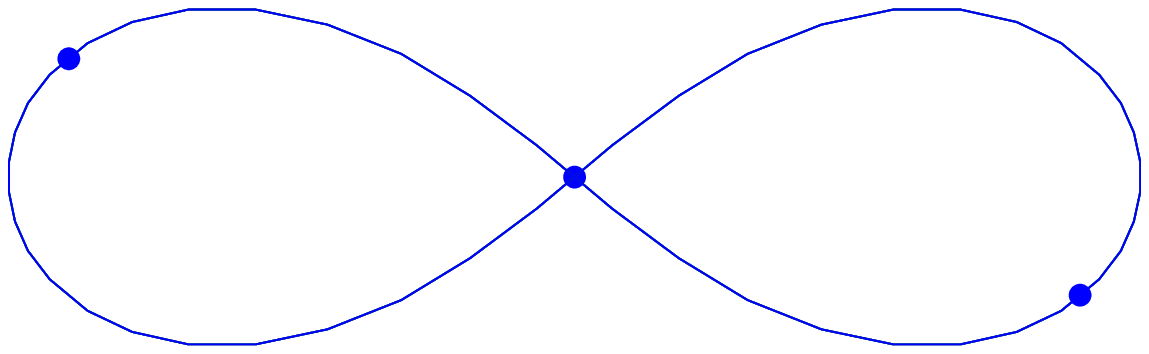}
    \caption{Figure-Eight Solution}\label{Eight}
    \end{minipage}
\end{figure}

Poincar\'e showed that it is impossible for $N\geq 3$ to find an
explicit expression for the general solution. So even for the
three-body problem, it is impossible to describe all the solutions.
Till now, for planar three-body problem, there are only three orbits
are found clearly for three equal masses, including the existence
and the shape of the orbits. They are Euler's solution (see
Fig.\ref{Euler}), Lagrange's solution (see Fig.\ref{Lagrange}) and
Figure-Eight solution (see Fig.\ref{Eight}). The Euler's solution
was discovered by Euler in 1765. It described a collinear
configuration, in which one body locates on the center of masses of
the three bodies, and the other two rotate around the center of the
masses on the two endpoints of the diameter of the circle orbit. The
Lagrange solution was discovered by Lagrange in 1772. The three
bodies form an equilateral triangle on a circle orbit and the center
of masses is the origin all the time. The Figure-Eight solution is a
new, surprisingly and remarkable solution after the two solutions,
which was found numerically by Moore first in 1993 (see
\cite{Moore}). In 2000, Chenciner and Montgomery rediscovered it and
proved its existence at the same time. It describes three equal
masses chasing one by one on an Eight-shaped curve in the plane. We
refer the reader to the following articles and references therein
for more detailed description: Chenciner et al \cite{CGMS},
Montgomery \cite{RM1} and Sim\'{o} \cite{{S1},{S2}} etc.

There is no other specific result for $N\geq 3$. Hence the
numerical solution is more and more important in finding and
describing new solutions vividly. In recent years, some new
periodic solutions of the N-body problem have been found by
minimizing the action of cycles with some special symmetric or
choreographic constraint (e.g. \cite{SU1},\cite{Zhang} etc), in
which Fourier approach is used by most authors for the numerical
solutions of N-body problem. In this paper, we propose the N-body
type difference equations, and from which we can discover some new
periodic orbits by minimizing the corresponding action functional.

The N-body type problems have been studied by Bahri, Rabinowitz
\cite{BR1}, Zelati \cite{VCZ1} among others. In \cite{VCZ1}, Zelati
proved the existence of periodic solutions under the radial symmetry
assumption. In this paper we prove the existence of periodic
solutions for the difference equations with multi-radial symmetric
constraint, which includes the radial symmetry. In the section on
numerical examples, we propose a new variational difference method
for the periodic solutions numerically.

The periodic solutions for N-body type problems (including N-body
problems) can be described by
\begin{equation}
\tag{N}\label{N} -m_{i}\ddot{q_{i}} =\nabla U(q),\quad i=1,\cdots,N,
\end{equation}
\begin{equation}
\tag{P}\label{P} q(0)=q(T), \quad \dot{q}(0)=\dot{q}(T),
\end{equation}
where $N$, $m_{i}$ and $q_{i}$ are defined as above,
$q=(q_{1},\cdots,q_{N})$ in $(R^{d})^{N}$, $T$ is the period, and
generally the singular potential function $U$ is
\begin{equation*}
U(q)=\frac{1}{2}\sum_{1\leq i\neq l\leq N}U_{il}(q_{l}-q_{i})
\end{equation*}
with assumptions: $\forall 1\leq i\neq l \leq N$ and $\forall \xi\in
R^{d}\backslash \{0\}$, we have
\begin{equation*}\tag{$A1$} \label{a1}
U_{il}\in C^{2}(R^{d}\backslash \{0\},R),\quad
U_{il}(\xi)=U_{ji}(\xi),
\end{equation*}

\begin{equation*}\tag{$A2$} \label{a2}
U_{il}(\xi) \rightarrow -\infty,\quad \mbox{as} \quad
|\xi|\rightarrow 0,
\end{equation*}

\begin{equation*}\tag{$A3$} \label{a3}
U(q_{1},\cdots,q_{N})\leq 0, \quad \forall(q_{1},\cdots,q_{N})\in
(R^{d})^{N},
\end{equation*}
where $U_{il}$ is usually described by
\begin{equation}\label{Uil1}
U_{il}=\frac{m_{i}m_{l}}{|q_{i}-q_{l}|^{\alpha}}.
\end{equation}

\section{Difference Equations}\label{sec:1}
In order to formulate the difference equations for N-body type
problems (\ref{N})-(\ref{P}), we discretize time within one period
$[0,T]$ with $k+1$ points $t_{0},\cdots,t_{k}$ and denote the set of
these discretized points as $S$, i.e. $ S=\{t_{0},\cdots,t_{k}\}$,
and introduce a difference scheme by
\begin{equation}
\tag{DS}\label{ds}
q(t)\doteq\hat{q}(t):=\frac{t_{s}-t}{h}q(s-1)+\frac{t-t_{s-1}}{h}q(s),
\quad s=1,\cdots,k+1,
\end{equation}
where $\hat{q}(s)$ is an abbreviation of $\hat{q}(t_{s})$, $t_{s}$
is the $s$th discretized point of $[0,T]$,
$\hat{q}=(\hat{q}_{1},\cdots,\hat{q}_{N})$, and $h$ is the
discretized step. Moreover, the traditional forth difference
scheme is used to approximate the derivative of $q$, which is used
in \cite{YJS02}. Then the first order difference of $q$ is
\begin{equation*}
\Delta q(t_{s})=q(s+1)-q(s),
\end{equation*}
and the second order difference of $q$ is
\begin{equation*}
\Delta^{2} q(t_{s})=q(s+1)-2q(s)+q(s-1).
\end{equation*}
Hence the derivatives of $q(t_{s})$ can be approximated by
\begin{equation*}
\dot{q}(t_{s})\doteq\frac{1}{h}\Delta
q(t_{s}),\quad\quad\ddot{q}(t_{s})\doteq\frac{1}{h^{2}}\Delta^{2}q(t_{s}),
\end{equation*}
By the difference approach above, (\ref{N}) can be transformed
into
\begin{equation*}
-\frac{1}{h^{2}}m_{i}(q_{i}(s+1)-2q_{i}(s)+q_{i}(s-1))=\nabla_{\hat{q}_{i}(s)}
\hat{U},
\end{equation*}
where $i=1,\cdots,N$ and
\begin{equation}\label{U}
\hat{U}=\sum_{s=1}^{k}U(\hat{q}(s))=\frac{1}{2}\sum_{s=1}^{k}\sum_{1\leq
i\neq l\leq N}U_{il}(\hat{q}_{l}(s)-\hat{q}_{i}(s)).
\end{equation}
For convenience, we set the following notations:
\begin{equation*}
M_{k} = \left( \begin{array}{cccccc} 2 & -1 & 0 & \cdots & -1 \\
-1 & 2& -1 & \ddots & \vdots \\
0 & \ddots & \ddots & \ddots & 0 \\
\vdots & \ddots & -1 & 2 & -1 \\
-1 & \cdots & 0 & -1 & 2 \end{array}\right),
\end{equation*}
\begin{equation*}
\hat{q} = \left( \begin{array}{c}
    \hat{q}_{1}^{\top}\\
    \hat{q}_{2}^{\top}\\
    \vdots \\
    \hat{q}_{N}^{\top}\\
    \end{array}
    \right),\quad
\hat{q}_{i}^{T} = \left( \begin{array}{c}
    \hat{q}_{i}^{(1)^{\top}}\\
    \hat{q}_{i}^{(2)^{\top}}\\
    \vdots \\
    \hat{q}_{i}^{(d)^{\top}}\\
    \end{array}
    \right),\quad
\hat{q}_{i}^{(p)^{T}} = \left( \begin{array}{c}
    \hat{q}_{i(1)}^{(p)^{\top}}\\
    \hat{q}_{i(2)}^{(p)^{\top}}\\
    \vdots \\
    \hat{q}_{i(k)}^{(p)^{\top}}\\
    \end{array}
    \right),
\end{equation*}
and
\begin{equation*}
M = \left( \begin{array}{ccc}
    M_{1} & & \\
     & \ddots & \\
    & & M_{N} \\
    \end{array}
    \right),
\end{equation*}
where
\begin{equation*}
M_{i} = \left( \begin{array}{ccc}
    m_{i}M_{k}&        &             \\
              & \ddots &             \\
              &        & m_{i}M_{k}  \\
    \end{array}
    \right),
\end{equation*}
\begin{equation*}
i=1,\cdots,N,\quad p=1,\cdots,d.
\end{equation*}
It is easy to verify that $M$ is a positive definitive matrix. Then
the difference equation of (\ref{N}) becomes
\begin{equation}
\tag{N1}\label{N1} \frac{1}{h}M \cdot
\hat{q}-\nabla_{\hat{q}_{i}(s)} \hat{U}h = 0,
\end{equation}
and it is obviously that the difference boundary conditions can be
described by
\begin{equation*}
\tag{P1}\label{P1} \hat{q}(j)=\hat{q}(k+j),\quad
j\in\mathbb{Z}^{+},
\end{equation*}
where $\mathbb{Z}^{+}$ is the set of positive integer. Hence the
Lagrangian functional of (\ref{N1})-(\ref{P1}) is the following
\begin{equation}\label{J}
J(q) = \frac{1}{2h}<M\hat{q},\hat{q}>-\hat{U}h.
\end{equation}

\begin{theorem}\label{p0}
There exists a small number $\delta>0$ for all $i\neq j$, $i, j=1,
\cdots, N$, $s\in S$, such that
$|\hat{q}_{i}(s)-\hat{q}_{j}(s)|>\delta$. Then the solution
$\hat{q}$ of (\ref{N1})-(\ref{P1}) approximates the solution $q$ of
(\ref{N})-(\ref{P}) when $k\rightarrow\infty$.
\end{theorem}
\begin{proof}
Since $|\hat{q}_{i}-\hat{q}_{j}|>\delta$, we can deduce the
following equation from (\ref{Uil1}) and (\ref{U}),
\begin{equation}\label{Ud}
\begin{array}{lll}
\nabla_{\hat{q}_{i}(s)} \hat{U}
&=&\displaystyle\frac{1}{2}\sum_{s=1}^{k}\sum_{1\leq
i\neq l\leq N}\cfrac{1}{|\hat{q}_{l}(s)-\hat{q}_{i}(s)|^{(\alpha+1)}}\\
&<&\displaystyle\frac{1}{2}\sum_{s=1}^{k}\sum_{1\leq
i\neq l\leq N}\frac{1}{\delta^{(\alpha+1)}}\\
&=&C,
\end{array}
\end{equation}
where $C$ is a constant. Similar to (\ref{ds}), we introduce the
difference scheme to $\dot{q}(t)$, i.e. we have
\begin{equation*}
\hat{\dot{q}}(t)=\frac{t_{s}-t}{h}\dot{q}(t_{s-1})+\frac{t-t_{s-1}}{h}\dot{q}(t_{s}),
\end{equation*}
where $t\in[t_{s-1}, t_{s}]$ and $t_{s-1},t_{s}\in S$. Notice that
\begin{equation*}
\dot{\hat{\dot{q}}}(t)=\frac{1}{h}(\dot{q}(t_{s})-\dot{q}(t_{s-1}))=\frac{1}{h}\Delta
\dot{q}(t_{s-1}).
\end{equation*}
Then when we let
\begin{equation*}
x=(q,\dot{q}),
\end{equation*}
the system (\ref{N})-(\ref{P}) can be rewritten as follows:
\begin{equation*}
\dot{x}=F(x),\quad x(0)=x(T),
\end{equation*}
where $x$ satisfies (\ref{ds}) and $F\in C(R^{2dN},R^{2dN})$, that
is
\begin{equation*}
F(x)=\left(
     \begin{array}{c}
     \dot{q}\\
     \ddot{q}
     \end{array}
     \right)=
     \left(
     \begin{array}{c}
     \vdots\\
     \dot{q}_{i}\\
     \vdots\\
     -\frac{1}{m_{i}}\nabla U(q)\\
     \vdots
     \end{array}
     \right).
\end{equation*}

From (\ref{Ud}), we have $|F|\leq M_{0}$ for some constant
$M_{0}$, and for $\forall\varepsilon>0$, there exist $\delta'>0$,
such that
\begin{equation*}
|F(x_{1})-F(x_{2})|\leq\varepsilon,
\end{equation*}
when $|x_{1}-x_{2}|\leq\delta'$. From (\ref{ds}), for
$t\in[t_{s-1},t_{s}]$, $s=1,\cdots,k$, clearly we have
\begin{eqnarray*}
\left|\hat{x}_{k}(t)-x(t_{s-1})\right|
&=&\left|\frac{t_{s}-t}{h}x(t_{s-1})+\frac{t-t_{s-1}}{h}x(t_{s})-x(t_{s-1})\right|\\
&=&\left|\frac{t-t_{s-1}}{h}\left(x(t_{s})-x(t_{s-1})\right)\right|\\
&=&\left|\frac{t-t_{s-1}}{h}F(x(t_{s-1}))(t_{s}-t_{s-1})\right|\\
&=&\left|(t-t_{s-1})F(x(t_{s-1}))\right|\\
&\leq&hM_{0}.
\end{eqnarray*}
Therefore, there exists a $K(\epsilon)>0$ such that
\begin{equation}\label{eq0}
|F(\hat{x}_{k}(t))-F(x(t_{s-1}))|\leq\varepsilon,~{\rm whenever}~
k\geq K(\epsilon).
\end{equation}

The error within one period can be estimated by
\begin{equation*}
\Delta_{k}(t)=\hat{x}_{k}(t)-\left(x(0)+\int_{0}^{t}F(\hat{x}_{k}(\tau))d\tau\right).
\end{equation*}
Hence
\begin{equation*}
\Delta_{k}(T)=\sum_{s=1}^{k}\int_{t_{s-1}}^{t_{s}}\left(F(x(t_{s-1}))
-F(\hat{x}_{k}(\tau))\right)d\tau.
\end{equation*}
From (\ref{eq0}), noting that the discretized step is $h=1$,  we
have
\begin{equation*}
|\Delta_{k}(T)|\leq\sum_{s=1}^{k}\varepsilon h =T\varepsilon.
\end{equation*}
This means $\Delta_{k}\rightarrow0$, $x_{k}\rightarrow x$, i.e.
$\hat{q}\rightarrow q$ when $k\rightarrow\infty$. The theorem is
proved.
\end{proof}

%\section{Difference Equation and Multi-Radial Symmetry}\label{sec:2}

Let
\begin{equation*}
\Omega=\{\hat{q}\ |\ \hat{q}\in(R^{d})^{N},\ \mbox{s.t.} \
\hat{q}_{i}\neq \hat{q}_{l},\ \forall i\neq l \}, \ \
\Lambda=\{\hat{q}\ |\ \hat{q}\in H^{1}(R/T;\Omega) \}, \ \
R^{d}=\prod_{i=1}^{b_{d}}R^{b_{i}}.
\end{equation*}
Then the {\em multi-radial symmetry} can be described as follows:
\begin{equation}
\tag{Ms}\label{Ms}\Lambda_{0}=\{\hat{q}\ |\ \hat{q}\in \Lambda,\
\hat{q}^{(b_{i})}(s+\tilde{k}^{i})=-\hat{q}^{(b_{i})}(s),
i=1,\cdots,b_{d}, \sum_{i=1}^{b_{d}}b_{i}=d \},
\end{equation}
where $\tilde{k}^{i}=A/A_{i}$, $A_{i}, A\in \mathbb{Z}^{+}$, $A$
is a common multiple of $A_{i}$, and
$\hat{q}^{(b_{i})}(s+\tilde{k}^{i})=-\hat{q}^{(b_{i})}(s)$ means
\begin{equation*}
\hat{q}^{(b_{i})}(t)=-\hat{q}^{(b_{i})}(\tau),\quad\mbox{ where }\ \
t\in[t_{s+\tilde{k}^{i}-1},t_{s+\tilde{k}^{i}}],\quad
\tau\in[t_{s-1},t_{s}].
\end{equation*}

\section{The Existence of Periodic Solutions}\label{sec:3}

Let $H^{1}$ is a Hilbert space, $\hat{q}\in
H^{1}([0,T],(R^{d})^{N})$, i.e. $\hat{q}_{i}\in
H^{1}([0,T],R^{d})$, $i=1,\cdots,N$, and we define the norm of
$\hat{q}$ in $H^{1}$ by
\begin{equation*}
\parallel\hat{q}\parallel=\sqrt{{S_{1}}^{2}+S_{2}},
\end{equation*}
where
\begin{equation*}
S_{1}=\frac{1}{T}\sum_{s=1}^{k}\int_{t_{s-1}}^{t_{s}}\left|\hat{q}(t)\right|dt,
\ \
S_{2}=\sum_{s=1}^{k}\int_{t_{s-1}}^{t_{s}}\left|\dot{\hat{q}}(t)\right|^{2}dt.
\end{equation*}

The following lemma describes the relationship between the
critical points of $J$ and the periodic solutions of
(\ref{N1})-(\ref{P1}).
\begin{lemma}\label{l1}
$(\hat{q}(1),\cdots,\hat{q}(k))^{\top}$ is a critical point of
$J$ if and only if $(\hat{q}(0),\cdots,\hat{q}(k+1))^{\top}$ is a
solution of (\ref{N1})-(\ref{P1}).
\end{lemma}
\begin{proof}
Since the difference equations transform the infinite dimensional
dynamics into finite dimensional dynamics, it is reasonable to view
$J$ as a continuously differentiable functional on $H^{1}$.

Let $\hat{q}=(\hat{q}(1),\cdots,\hat{q}(k))^{\top}$ is the critical
point of $J$ if and only if the first order derivative of $J$ on
$\hat{q}$ is equal to zero. That is
\begin{equation*}
J^{'}(\hat{q})=M\hat{q}-\nabla_{\hat{q}}\hat{U}=0.
\end{equation*}

Notice that
\begin{equation*}
(\hat{q}(0),\hat{q}(1),\cdots,\hat{q}(k),\hat{q}(k+1))^{\top}=
(\hat{q}(k),\hat{q}(1),\cdots,\hat{q}(k),\hat{q}(1))^{\top}.
\end{equation*}
This satisfies (\ref{N1})-(\ref{P1}) obviously. Thus the proof is
completed.
\end{proof}

\begin{remark}
Actually, $\hat{q}_{i}$ is a $k$-polygon in $R^{d}$ as a simple case
with the difference scheme (\ref{ds}). Hence we can call
$\{\hat{q}(s)\}_{s=0}^{k+1}=(\hat{q}(0),\cdots,\hat{q}(k+1))^{T}$ a
`cycle' in $(R^{d})^{N}$, which means obviously that
$\{\hat{q}_{i}(s)\}_{s=0}^{k+1}=(\hat{q}_{i}(0),\cdots,\hat{q}_{i}(k+1))^{T}$,
 $i=1,\cdots,N$, are `cycles' in $R^{d}$.
\end{remark}

We denote
\begin{equation*}
\Delta =\bigcup_{s=1}^{k}\{t\in[t_{s-1},t_{s}] \ | \
\hat{q}_{i}(t)=\hat{q}_{l}(t),\ \mbox{for some}\ i \neq l \}.
\end{equation*}

\begin{definition}We call $\hat{q}$ is a {\em non-collision solution} of (\ref{N1}) and (\ref{P1}),
if $\Delta=\phi$.
\end{definition}
\begin{definition}We call $\hat{q}$ is a {\em generalized solution} of (\ref{N1}) and (\ref{P1}), if

a) $meas(\Delta)=0$;

b) $\hat{q}$ satisfy (\ref{N1});

c)
\begin{equation}\label{def2}
 \sum_{i=1}^{N}\frac{m_{i}}{2}\left|\dot{\hat{q}}_{i}(t)\right|^{2}-U(\hat{q}(t))\equiv
C,
\end{equation}
where $t\in[0,T]\backslash \Delta$ and $C$ is a constant.
\end{definition}
\begin{remark}
It is easy to verify that the integral of (\ref{def2}) with respect
to $t$ approximates $J(\hat{q})$ since
\begin{equation}
\int_{t_{s-1}}^{t_{s}}U_{il}(\hat{q}_{l}(t)-\hat{q}_{i}(t))dt\doteq
U_{il}(\hat{q}_{l}(s)-\hat{q}_{i}(s)).
\end{equation}
\end{remark}

For convenience, we suppose $\hat{U}$ in (\ref{N1}) has only one
singular point, and without loss of generality, let it be $0\in
R^{d}$. Hence we have

\begin{definition} Let some $V_{il}\in C^{1}(R^{d} \backslash \{0 \} , R )$,
$V_{il}\rightarrow +\infty$ ($\xi \rightarrow 0$) and
\begin{equation}
\tag{SF} \label{sf} -U_{il}(\xi)\geq |\nabla V_{il}(\xi)|^{2},
\qquad \forall \xi \in R^{d} \backslash \{0 \},
\end{equation}
where $|\xi|$ is small. Then we say $U_{il}$ satisfies {\em Strong
Force {\rm (SF)} condition}.
\end{definition}

The condition (\ref{sf}) is very efficient to overcome the
obstacle coming from singular potential. In \cite{WG}, (\ref{sf})
is used to bound $q$ away from $S$, where $S$ is a closed nonempty
set, and $U\rightarrow -\infty$ as $q\rightarrow S$. Here, we
introduce Lemma 2 based on the discussions above for N-body type
difference equations.

Let $\Gamma$ be the family of homotopic `cycles' in $H^{1}([0,T],
(R^{d}\backslash\{0\})^{N})$, then we have

\begin{lemma}\label{l2}
Assume condition (\ref{sf}) is satisfied, and there exists a
constants $c$ such that
\begin{equation*}
\langle M\hat{q},\hat{q}\rangle\leq c, \quad|\hat{U}|\leq c,
\end{equation*}
for all $\hat{q}\in \Gamma$, then $\hat{q}$ is bounded away from
$\{0\}$. That is, there exists $\delta>0$ such that no cycles in
$\Gamma$ intersect the $\delta$-neighborhood of $\{0\}$.
\end{lemma}

\begin{proof}
From (\ref{sf}), it is easy to verify that $U_{il}$ are all
negative, then from (\ref{U}) we have the following estimate for
$\hat{U}$:
\begin{equation*}
|\hat{U}|=\frac{1}{2}\sum_{s=1}^{k}\sum_{1\leq i\neq l\leq
N}|U_{il}(\hat{q}_{l}(s)-\hat{q}_{i}(s))|\leq c.
\end{equation*}

Let
\begin{equation*}
B = \left( \begin{array}{ccccccc}
m_{1}  &        &       &        &       &        &      \\
       & \ddots &       &        &       &        &      \\
       &        & m_{1} &        &       &        &      \\
       &        &       & \ddots &       &        &      \\
       &        &       &        & m_{N} &        &      \\
       &        &       &        &       & \ddots &      \\
       &        &       &        &       &        & m_{N}\\
\end{array}\right),
\end{equation*}
\begin{equation*}
\tilde{m}=\min\{m_{1},\cdots,m_{N}\},\quad y=B\hat{q},
\end{equation*}
and let
\begin{equation*}
G=\{\rho=(y,z)\in(R^{d})^{N}\times R\ | \ \
y,\hat{q}\in(R^{d}\backslash\{0\})^{N},\ \ z=V(\hat{q})\}
\end{equation*}
be the graph of $V$ in (\ref{sf}). Suppose the binary relation is
\begin{equation*}
\rho=(y,V(\hat{q})),
\end{equation*}
whose arc length is $a=a(\rho)$. Then we have
\begin{equation*}
\begin{array}{lll}
a
&=& \int_{0}^{T}|\dot{a}|dt=\sum\limits_{s=1}^{k}\int_{t_{s-1}}^{t_{s}}|\dot{a}|dt  \\
&=&\sum\limits_{s=1}^{k}\int_{t_{s-1}}^{t_{s}}\sqrt{|\dot{y}|^{2}+\langle\hat{q},\nabla V\rangle^{2}}dt  \\
&\leq&
\sum\limits_{s=1}^{k}\int_{t_{s-1}}^{t_{s}}\left(|\dot{y}|^{2}+|\dot{\hat{q}}|^{2}|\nabla
V|^{2}\right)^{\frac{1}{2}}dt \\
&\leq&
\sum\limits_{s=1}^{k}\int_{t_{s-1}}^{t_{s}}\left(|\dot{y}|^{2}\right)^{\frac{1}{2}}\left(1+\frac{1}{\tilde{m}^{2}}|\nabla
V|^{2}\right)^{\frac{1}{2}}dt \\
&\leq&
\sum\limits_{s=1}^{k}\left(\int_{t_{s-1}}^{t_{s}}|\dot{y}|^{2}dt\right)^{\frac{1}{2}}\left(\int_{t_{s-1}}^{t_{s}}\left(1+\frac{1}{\tilde{m}^{2}}|\nabla
V|^{2}\right)dt\right)^{\frac{1}{2}} \quad \mbox{from H\"{o}lder inequality} \\
&\leq&
\left(\sum\limits_{s=1}^{k}\int_{t_{s-1}}^{t_{s}}|\dot{y}|^{2}dt\right)^{\frac{1}{2}}\left(\sum\limits_{s=1}^{k}\int_{t_{s-1}}^{t_{s}}\left(1+\frac{1}{\tilde{m}^{2}}|\nabla
V|^{2}\right)dt\right)^{\frac{1}{2}} \\
&=&
\left(\sum\limits_{s=1}^{k}\int_{t_{s-1}}^{t_{s}}|B\dot{\hat{q}}|^{2}dt\right)^{\frac{1}{2}}\left(T+\frac{1}{\tilde{m}^{2}}\sum\limits_{s=1}^{k}\int_{t_{s-1}}^{t_{s}}\left(|\nabla
V|^{2}\right)dt\right)^{\frac{1}{2}} \\
&\leq&
\left(\sum\limits_{s=1}^{k}\int_{t_{s-1}}^{t_{s}}|B\dot{\hat{q}}|^{2}dt\right)^{\frac{1}{2}}\left(T+\frac{1}{\tilde{m}^{2}}\sum\limits_{s=1}^{k}\int_{t_{s-1}}^{t_{s}}|U|dt\right)^{\frac{1}{2}} \quad \mbox{ from (\ref{sf}) condition} \\
&=&
\left(\sum\limits_{s=1}^{k}|B(\hat{q}(s)-\hat{q}(s-1))|^{2}\right)^{\frac{1}{2}}\left(T+\frac{1}{\tilde{m}^{2}}\sum\limits_{s=1}^{k}\int_{t_{s-1}}^{t_{s}}|U(\hat{q}(t))|dt\right)^{\frac{1}{2}} \\
&=&
\left(\sum\limits_{s=1}^{k}\sum\limits_{i=1}^{N}m_{i}|(\hat{q}_{i}(s)-\hat{q}_{i}(s-1))|^{2}\right)^{\frac{1}{2}} \\
&&\left(T+\frac{1}{\tilde{m}^{2}}\frac{1}{2}\sum\limits_{s=1}^{k}\sum\limits_{1\leq i\neq l\leq N}|U_{il}(\hat{q}_{l}(s)-\hat{q}_{i}(s))|\right)^{\frac{1}{2}} \\
&=&\sqrt{\langle M\hat{q},\hat{q}\rangle}\sqrt{T+\frac{1}{\tilde{m}^{2}}|\hat{U}|}\\
&\leq&\sqrt{c}\sqrt{T+\frac{1}{\tilde{m}^{2}}c}\equiv C \mbox{ (a
constant)}.
\end{array}
\end{equation*}

Now let $\{\hat{q}^{\nu}\}_{\nu=1}^{\infty}$ be a sequence in
$\Gamma$. If $\hat{q}^{\nu}$ are not bounded away from $\{0\}$, for
$\forall\varepsilon>0$, the $\varepsilon$-neighborhood of $\{0\}$ in
$R^{d}$, $B_{\varepsilon}$, must intersect with infinite many
$\hat{q}^{\nu}$. By passing to a subsequence, we can suppose that
\begin{equation*}
\inf|\hat{q}^{\nu}|\rightarrow 0,\quad\mbox{as}\quad
n\rightarrow\infty.
\end{equation*}
This implies that the corresponding $\rho^{\nu}$ ascend infinitely
far up to the skylight of $G$ at $\{0\}$ as $n\rightarrow \infty$.

Suppose there exists $\delta>0$, $B_{\delta}$ is
$\delta$-neighborhood of $\{0\}$, such that portions of
$\hat{q}^{\nu}$ fall outside the ball $B_{\delta}$ and the others
(infinite many) $\hat{q}^{\nu}$ in $B_{\delta}$. Then the
variation of $V(\hat{q}^{\nu})$ would become infinite as
$n\rightarrow \infty$, but $a(\rho)$ is greater than the variation
of $V(\hat{q}^{\nu})$, a contradiction to the fact $a\leq C$.

On the other hand, if there is no such $\delta$, $\hat{q}^{\nu}$
will collapse to the point $\{0\}$, which implies that
\begin{equation*}
\hat{U}\rightarrow +\infty.
\end{equation*}
This contradicts to the hypotheses at the beginning. Hence $\hat{U}$
is bounded away from $\{0\}$.

Now, suppose $\hat{q}^{\nu}$ converge uniformly to $\hat{q}^{*}$
in the norm in $H^{1}$. By the above argument, if $\hat{q}^{\nu}$
are not bounded away from $\{0\}$, $\hat{q}^{*}$ will be attached
to $\{0\}$. The proof is complete.
\end{proof}

\begin{remark}\ \

(a) For the case that the potential $U$ has more than one singular
points has the same result in Lemma 2.

(b) Lemma 2 means that there is no $\hat{q}\in \Gamma$ intersects
the $\delta$-neighborhood of $\{0\}$ under (\ref{sf}) condition.

(c) The converse-negative proposition of Lemma 2 shows that if
$\bar{q}\in\partial\Lambda$, $q^{\mu}\rightarrow\bar{q}$ weakly in
$H^{1}$ and strongly in $C^{0}$, then $J(q^{\mu})\rightarrow
+\infty$.
\end{remark}

With respect to the periodic solutions of (\ref{N1})-(\ref{P1}), we
have the following theorem:

\begin{theorem}\label{main}
Let $U_{il}$ and $\hat{U}$ in (\ref{U}) satisfy
(\ref{a1})-(\ref{a3}) and $\hat{q}$ satisfy (\ref{Ms}). Then
(\ref{N1})-(\ref{P1}) has infinitely many generalized solutions;
furthermore, if (\ref{sf}) holds, then (\ref{N1})-(\ref{P1}) has
infinitely many non-collision solutions.
\end{theorem}

\begin{proof}
It is obvious that $\Lambda_{0}$ is a subset of $\Lambda$, and
$\tilde{U}$ is an even function on $\Lambda_{0}$ for
$U_{il}=U_{li}$, so we can easily check that the critical points
of (\ref{Jd}) on $\Lambda_{0}$ are also the critical points on
$\Lambda$.

Since $U_{il}$ does not satisfy (\ref{sf}), we can modify it as
\begin{equation*}
U^{\delta}_{il}=U_{il}-\frac{\phi(|\hat{q}_{l}-\hat{q}_{i}|)}{|\hat{q}_{l}-\hat{q}_{i}|^{2}},
\end{equation*}
where $\delta$ is small and
\begin{displaymath}
\phi = \left\{
    \begin{array}{ll}
    0, & \quad\mbox{when} \ |\hat{q}_{l}-\hat{q}_{i}|\leq \displaystyle\frac{\delta}{2},\\
    1, & \quad\mbox{when} \ |\hat{q}_{l}-\hat{q}_{i}|\geq \delta.
    \end{array} \right.
\end{displaymath}
Then we get a new Lagrangian functional with condition (\ref{sf})
from (\ref{J}),
\begin{equation}\label{Jd}
J^{\delta}(\hat{q}) = \frac{1}{2} <M\hat{q},\hat{q}
>-\tilde{U}^{\delta},
\end{equation}
where
\begin{equation*}
\tilde{U}^{\delta} = \sum_{s=1}^{k} \sum_{1\leq i < l \leq N}
U_{il}^{\delta}.
\end{equation*}

Let $J^{\delta}$ have infimum
\begin{equation*}
c_{\delta}=\inf \{ \ J^{\delta}(\hat{q}) \ | \ \hat{q}\in
\Lambda_{0} \ \}.
\end{equation*}

Consider minimizing sequence $\hat{q}^{\mu}(s)\in \Lambda_{0}$ such
that $J^{\delta}(\hat{q}^{\mu}(s))\rightarrow c_{\delta}$ when
$\mu\rightarrow +\infty$. Then, when $\mu$ is large enough, we have
\begin{equation*}\label{Jd}
J^{\delta}(\hat{q}^{\mu}(s))\leq c_{\delta}+1.
\end{equation*}

Since $\hat{q}\in\Lambda_{0}$, we have that
\begin{eqnarray*}
\left|\hat{q}_{i}(s)\right| & \leq &
d\times\max\{|\hat{q}^{b}_{i}(s)|,|\hat{q}^{b'}_{i}(s)|\} \\
& =  &
d\times\max\{\frac{1}{2}|\hat{q}^{b}_{i}(s)-\hat{q}^{b}_{i}(s+\tilde{k})|,
\frac{1}{2}|\hat{q}^{b'}_{i}(s)-\hat{q}^{b'}_{i}(s+\hat{k})|\} \\
&\leq&
d\times\max\{\frac{1}{2}\sum_{j=0}^{\tilde{k}-1}|\hat{q}_{i}(s-j)-\hat{q}_{i}(s-j+1)|,\\
&&\quad\quad\quad\quad\frac{1}{2}\sum_{j=0}^{\hat{k}-1}|\hat{q}_{i}(s-j)-\hat{q}_{i}(s-j+1)|\}  \\
&\leq&
\frac{d}{2}\sum_{j=0}^{k-1}|\hat{q}_{i}(s-j)-\hat{q}_{i}(s-j+1)|\\
&=&
\frac{d}{2}\sum_{j=0}^{k-1}\left|\hat{q}_{i}(s-j)-\hat{q}_{i}(s-j+1)\right|  \\
\end{eqnarray*}
\begin{eqnarray*}
& =  &
\frac{d}{2}\sum_{j=0}^{k-1}\left(\sqrt{m_{i}}\left|\hat{q}_{i}(s-j)-\hat{q}_{i}(s-j+1)\right|\sqrt{\frac{1}{m_{i}}}\right)\\
&\leq&
\frac{d}{2}\left(\sum_{j=0}^{k-1}m_{i}\left|\hat{q}_{i}(s-j)-\hat{q}_{i}(s-j+1)\right|^{2}\right)^{\frac{1}{2}}
\left(\sum_{j=0}^{k-1}\frac{1}{\tilde{m}}\right)^{\frac{1}{2}}\\
&=&\frac{d}{2}\left(<M\hat{q},\hat{q}>\right)^{\frac{1}{2}}\left(\frac{T}{\tilde{m}}\right)^{\frac{1}{2}}\\
&\leq&d\times\frac{1}{2}\left(J^{\delta}\right)^{\frac{1}{2}}\left(\frac{T}{\tilde{m}}\right)^{\frac{1}{2}} \\
&\leq&
\frac{d}{2}\left(\frac{T(c_{\delta}+1)}{\tilde{m}}\right)^{\frac{1}{2}},
\quad s=1,\cdots,k.
\end{eqnarray*}
where $\tilde{m}=\min\limits_{i=1,\cdots,N} \{ m_{i}\}$. This
illuminates that there exists a $\bar{q}\in \Lambda_{0}$ and a
subsequence $\{\hat{q}^{\mu}\}_{\mu}$ such that
\begin{equation*}
\hat{q}^{\mu}\rightarrow\bar{q}(s),\quad\quad s=1,\cdots,k,
\end{equation*}
\begin{equation*}
J^{\delta}(\bar{q})=c_{\delta}.
\end{equation*}
Then $\bar{q}$ is a critical point of (\ref{Jd}).

Now, from the remarks of Lemma \ref{l2} and Lemma \ref{l1}, there
exist a subsequence $\{\hat{q}^{\mu}(s)\}$ and
$\hat{q}^{*}_{\delta}$ in $H^{1}(R/T,(R^{d})^{N})$, such that
$\hat{q}^{*}_{\delta}$ is a non-collision periodic solution of
N-body problem (\ref{N1})-(\ref{P1}).

Along the line of Theorem1.1 in \cite{VCZ1}, we can find a
generalized solution $\hat{q}^{*}$ of (\ref{N1})-(\ref{P1}) by
constructing compact nested interval cover of $[0,T]\backslash
\Delta$ to show $U^{\delta}\rightarrow U$ and
$\hat{q}^{*}_{\delta}\rightarrow \hat{q}^{*}$ when
$\delta\rightarrow 0$. We can also prove $Q_{T/C}$ is a
$T$-periodic solution of (\ref{N1})-(\ref{P1}) if $Q$ is, where
$C$ is a positive integer, then we can find infinitely many
$T$-periodic solutions of (\ref{N1})-(\ref{P1}).  This completes
the proof.
\end{proof}

\section{Numerical Examples for the N-Body Problem}\label{sec:4}
In this section, we present some numerical solutions of N-body
problem (\ref{Nb})-(\ref{P}) as an example of N-body type. In this
case, $\alpha=1$ in (\ref{Uil1}), which is usually called
Newtonian Potential.

The method of steepest descent is adopted in finding the minimum
of the functional $J$ since this method is simple, easy to apply,
and each iteration is fast. If the minimum points exist, the
method guarantees that we can locate them after at least an
infinite number of iterations. We implement it in MATLAB7.0 for
getting the following orbits. And from the process in deducing
$J$, we named the method \textit{variational difference method}.
The following examples show that it is at least very valid for our
problem.

When we set the multi-radial symmetry with
\begin{equation*}
q^{(1,2)}\left(t+\frac{T}{2}\right)=-q(t)^{(1,2)},\quad q^{(1,2)}\in
R^{2},
\end{equation*}
which is a simple radial symmetric constraint in the plane. We get
the following two figures with different initial values for
three-body problems with equal masses.

In the following tables, $\hat{q}_{i}$ in the columns is the initial
value of $i$th body, where $i=1,2,3$. And the $j$th raw is
corresponding to the initial value of the $j$th difference point,
where $j=1,\cdots,k$. We set the values at random by MATLAB.

\begin{table}[H]
\begin{minipage}[t]{0.5\linewidth}
\centering
\begin{tabular}{|c|c|c|}
\hline
$\hat{q}_{1}$ & $\hat{q}_{2}$ & $\hat{q}_{3}$ \\
\hline
$\begin{array}{c}4.8375\\4.3081\end{array}$&$\begin{array}{c}-3.5004\\-5.8243\end{array}$&$\begin{array}{c}-1.3367\\1.5174\end{array}$\\
\hline
$\begin{array}{c}5.0723\\3.9921\end{array}$&$\begin{array}{c}-3.3872\\-5.7708\end{array}$&$\begin{array}{c}-1.6835\\1.7809\end{array}$\\
\hline
\vdots & \vdots & \vdots\\
\hline
$\begin{array}{c}-4.2971\\-4.8963\end{array}$&$\begin{array}{c}3.7003\\5.8687\end{array}$&$\begin{array}{c}0.5988\\-0.9714\end{array}$\\
\hline
$\begin{array}{c}-4.5792\\-4.6106\end{array}$&$\begin{array}{c}3.6051\\5.8569\end{array}$&$\begin{array}{c}0.9749\\-1.2464\end{array}$\\
\hline
\end{tabular}
\parbox{2in}{\caption{The initial value of
Fig.\ref{L_100_0126}}\label{t2}}
\end{minipage}
\begin{minipage}[t]{0.5\linewidth}
\centering
\begin{tabular}{|c|c|c|}
\hline
$\hat{q}_{1}$ & $\hat{q}_{2}$ & $\hat{q}_{3}$ \\
\hline
$\begin{array}{c}4.8375\\4.3081\end{array}$&$\begin{array}{c}-3.5004\\-5.8243\end{array}$&$\begin{array}{c}-1.3367\\1.5174\end{array}$\\
\hline
$\begin{array}{c}5.0723\\3.9921\end{array}$&$\begin{array}{c}-3.3872\\-5.7708\end{array}$&$\begin{array}{c}-1.6835\\1.7809\end{array}$\\
\hline
\vdots & \vdots & \vdots\\
\hline
$\begin{array}{c}-4.2971\\-4.8963\end{array}$&$\begin{array}{c}3.7003\\5.8687\end{array}$&$\begin{array}{c}0.5988\\-0.9714\end{array}$\\
\hline
$\begin{array}{c}-4.5792\\-4.6106\end{array}$&$\begin{array}{c}3.6051\\5.8569\end{array}$&$\begin{array}{c}0.9749\\-1.2464\end{array}$\\
\hline
\end{tabular}
\parbox{2in}{\caption{The initial value of
Fig.\ref{N_104_0126}}\label{t3}}
\end{minipage}
\end{table}

\begin{figure}[H]
    \begin{minipage}[t]{0.5\linewidth}
    \centering
    \includegraphics[width=2.7in, height=2.6in]{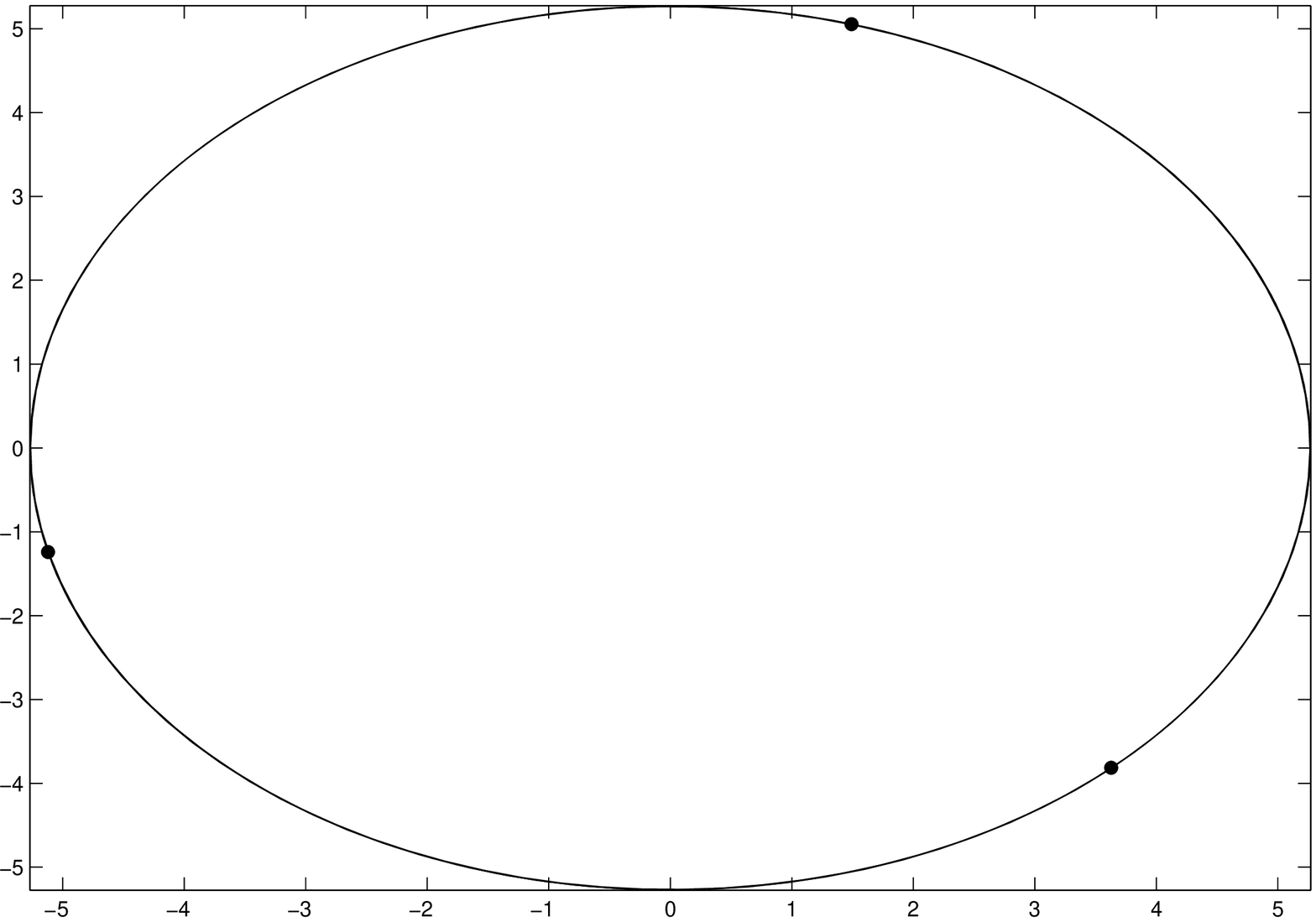}
    \parbox{2in}{\caption{$m=[1,1,1]$, $J=0.493$}\label{L_100_0126}}
    \end{minipage}
    \begin{minipage}[t]{0.5\linewidth}
    \centering
    \includegraphics[width=2.7in, height=2.6in]{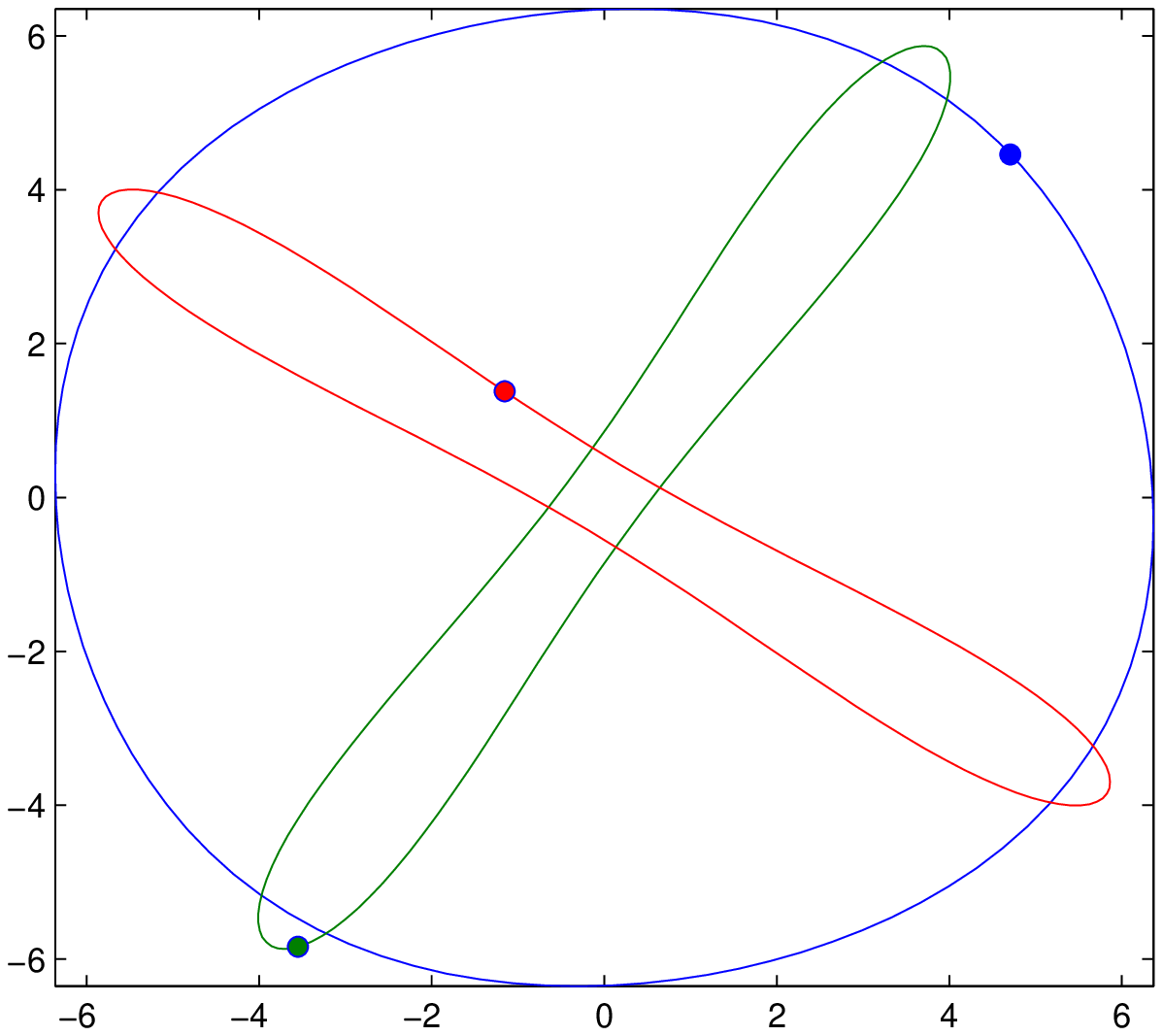}
    \parbox{2in}{\caption{$m=[1,1,1]$, $J=0.512$}\label{N_104_0126}}
    \end{minipage}
\end{figure}

When the distribution of the masses are
$m_{1}=1,m_{2}=0.5,m_{3}=1.5$, four different numerical solutions
are found with different initial values of $J$ for three-body
problem.
\begin{table}[H]
\begin{minipage}[t]{0.5\linewidth}
\centering
\begin{tabular}{|c|c|c|}
\hline
$\hat{q}_{1}$ & $\hat{q}_{2}$ & $\hat{q}_{3}$ \\
\hline
$\begin{array}{c}0.9977\\0.3897\end{array}$&$\begin{array}{c}0.7023\\0.5013\end{array}$&$\begin{array}{c}0.9878\\0.5027\end{array}$\\
\hline
$\begin{array}{c}0.9878\\0.5027\end{array}$&$\begin{array}{c}0.3214\\0.2917\end{array}$&$\begin{array}{c}0.8118\\0.2859\end{array}$\\
\hline
\vdots & \vdots & \vdots\\
\hline
$\begin{array}{c}0.8906\\0.4815\end{array}$&$\begin{array}{c}0.0109\\0.0134\end{array}$&$\begin{array}{c}0.1405\\0.1421\end{array}$\\
\hline
$\begin{array}{c}0.5134\\0.4002\end{array}$&$\begin{array}{c}0.0484\\0.0046\end{array}$&$\begin{array}{c}0.1205\\0.3442\end{array}$\\
\hline
\end{tabular}
\parbox{2in}{\caption{The initial value of Fig. \ref{P3_N_5_046_0207}}\label{t4}}
\end{minipage}
\begin{minipage}[t]{0.5\linewidth}
\centering
\begin{tabular}{|c|c|c|}
\hline
$\hat{q}_{1}$ & $\hat{q}_{2}$ & $\hat{q}_{3}$ \\
\hline
$\begin{array}{c}0.2916\\0.3066\end{array}$&$\begin{array}{c}0.2886\\0.3100\end{array}$&$\begin{array}{c}0.5938\\0.3580\end{array}$\\
\hline
$\begin{array}{c}0.0974\\0.7207\end{array}$&$\begin{array}{c}0.9711\\0.1348\end{array}$&$\begin{array}{c}0.1567\\0.9382\end{array}$\\
\hline
\vdots & \vdots & \vdots\\
\hline
$\begin{array}{c}0.6674\\0.1623\end{array}$&$\begin{array}{c}0.6226\\0.6986\end{array}$&$\begin{array}{c}0.3063\\0.3450\end{array}$\\
\hline
$\begin{array}{c}0.3109\\0.0311\end{array}$&$\begin{array}{c}0.1326\\0.8893\end{array}$&$\begin{array}{c}0.6609\\0.8840\end{array}$\\
\hline
\end{tabular}
\parbox{2in}{\caption{The initial value of
Fig. \ref{P3_N_1_030_0207}}\label{t5}}
\end{minipage}
\end{table}

\begin{figure}[H]
    \begin{minipage}[t]{0.5\linewidth}
    \centering
    \includegraphics[width=2.7in]{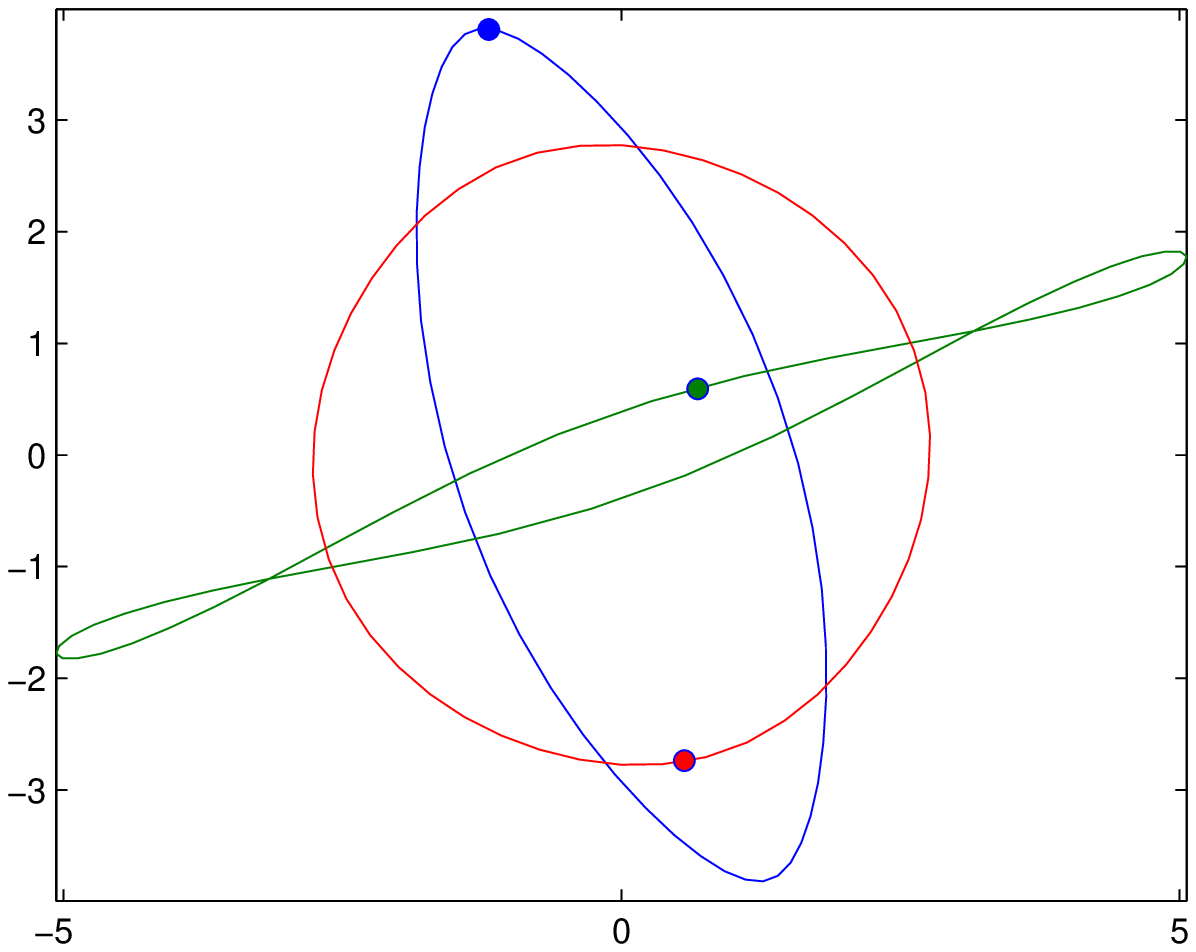}
    \parbox{2in}{\caption{$m=[1,0.5,1.5]$, $J=0.805$}\label{P3_N_5_046_0207}}
    \end{minipage}
%    \hspace{0.5cm}
    \begin{minipage}[t]{0.5\linewidth}
    \centering
    \includegraphics[width=2.7in]{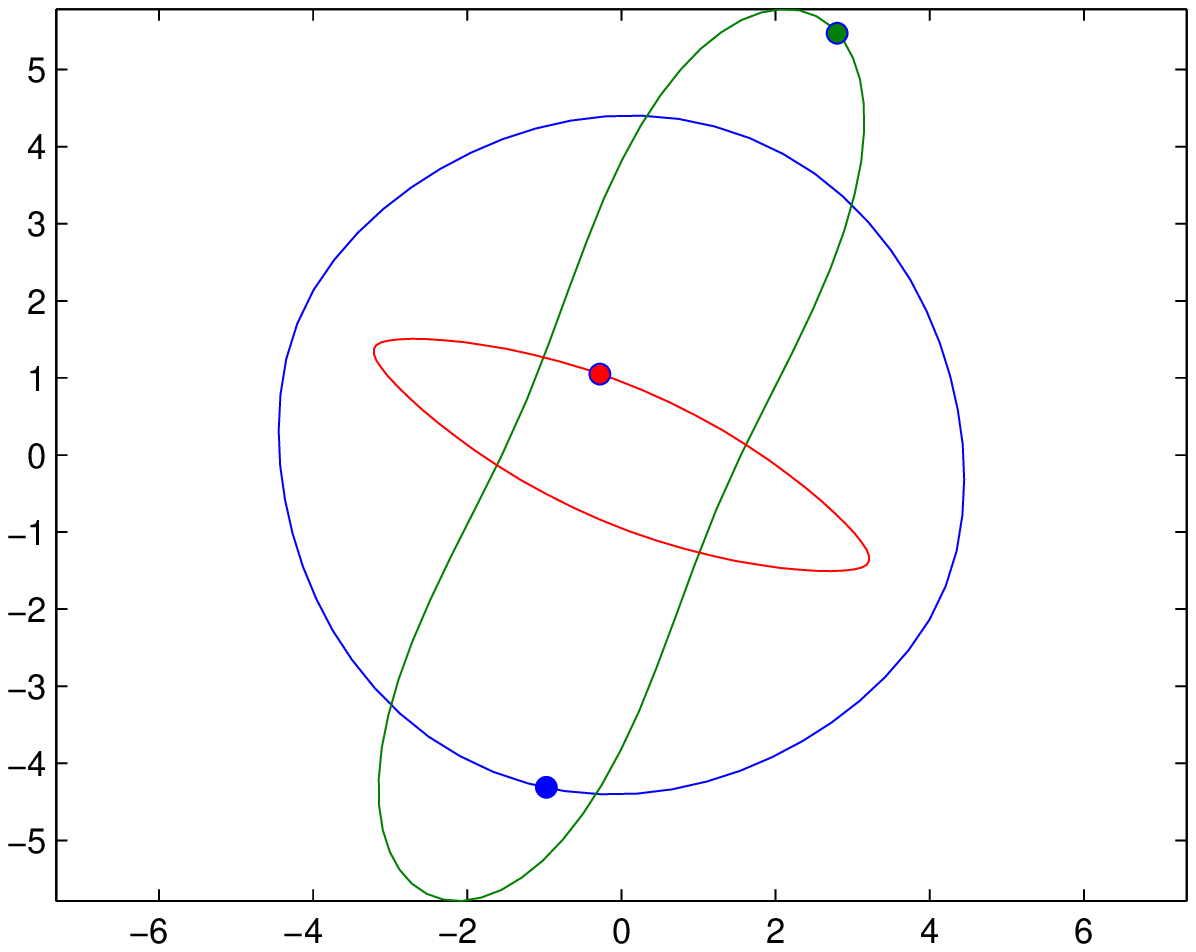}
    \parbox{2in}{\caption{$m=[1,0.5,1.5]$, $J=0.670$}\label{P3_N_1_030_0207}}
    \end{minipage}
\end{figure}

\begin{table}[H]
\begin{minipage}[t]{0.5\linewidth}
\centering
\begin{tabular}{|c|c|c|}
\hline
$\hat{q}_{1}$ & $\hat{q}_{2}$ & $\hat{q}_{3}$ \\
\hline
$\begin{array}{c}0.5499\\0.6606\end{array}$&$\begin{array}{c}0.8989\\0.8988\end{array}$&$\begin{array}{c}0.7706\\0.2469\end{array}$\\
\hline
$\begin{array}{c}0.0424\\0.2721\end{array}$&$\begin{array}{c}0.0684\\0.7092\end{array}$&$\begin{array}{c}0.6264\\0.1738\end{array}$\\
\hline
\vdots & \vdots & \vdots\\
\hline
$\begin{array}{c}0.7811\\0.7249\end{array}$&$\begin{array}{c}0.2253\\0.9763\end{array}$&$\begin{array}{c}0.1856\\0.1436\end{array}$\\
\hline
$\begin{array}{c}0.8518\\0.3835\end{array}$&$\begin{array}{c}0.2435\\0.5769\end{array}$&$\begin{array}{c}0.2441\\0.7770\end{array}$\\
\hline
\end{tabular}
\parbox{2in}{\caption{The initial value of
Fig. \ref{P3_N_2_046_0207}}\label{t6}}
\end{minipage}
\begin{minipage}[t]{0.5\linewidth}
\centering
\begin{tabular}{|c|c|c|}
\hline
$\hat{q}_{1}$ & $\hat{q}_{2}$ & $\hat{q}_{3}$ \\
\hline
$\begin{array}{c}0.9342\\0.4225\end{array}$&$\begin{array}{c}0.3759\\0.6273\end{array}$&$\begin{array}{c}0.1146\\0.4319\end{array}$\\
\hline
$\begin{array}{c}0.2644\\0.8560\end{array}$&$\begin{array}{c}0.0099\\0.6991\end{array}$&$\begin{array}{c}0.6649\\0.6343\end{array}$\\
\hline
\vdots & \vdots & \vdots\\
\hline
$\begin{array}{c}0.1556\\0.6252\end{array}$&$\begin{array}{c}0.6555\\0.7036\end{array}$&$\begin{array}{c}0.7505\\0.6802\end{array}$\\
\hline
$\begin{array}{c}0.1911\\0.7334\end{array}$&$\begin{array}{c}0.3919\\0.4850\end{array}$&$\begin{array}{c}0.7400\\0.0534\end{array}$\\
\hline
\end{tabular}
\parbox{2in}{\caption{The initial value of Fig. \ref{P3_N_1_040_0208}}\label{t7}}
\end{minipage}
\end{table}
\begin{figure}[H]
    \begin{minipage}[t]{0.5\linewidth}
    \centering
    \includegraphics[width=2.7in]{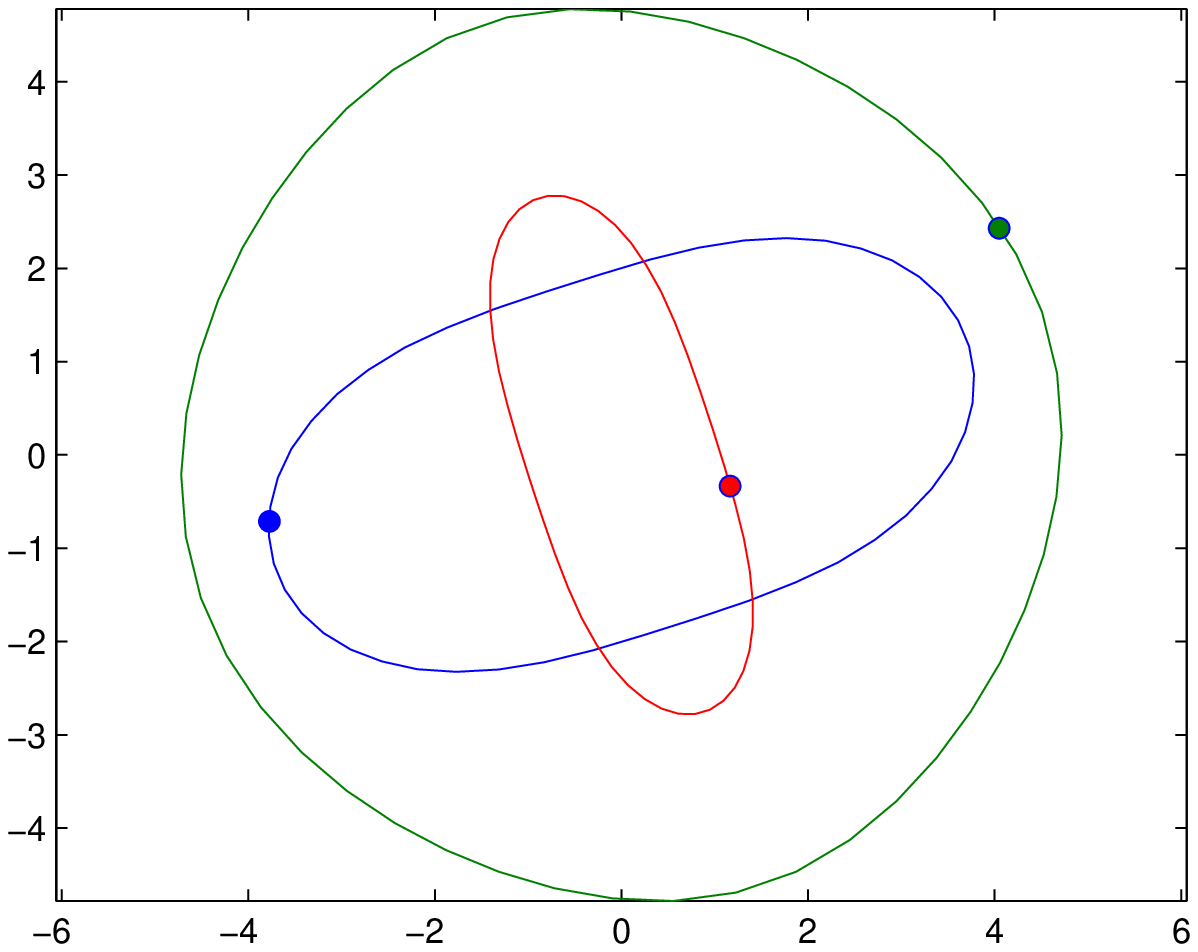}
    \parbox{2in}{\caption{$m=[1,0.5,1.5]$, $J=0.795$}\label{P3_N_2_046_0207}}
    \end{minipage}
    \begin{minipage}[t]{0.5\linewidth}
    \centering
    \includegraphics[width=2.7in]{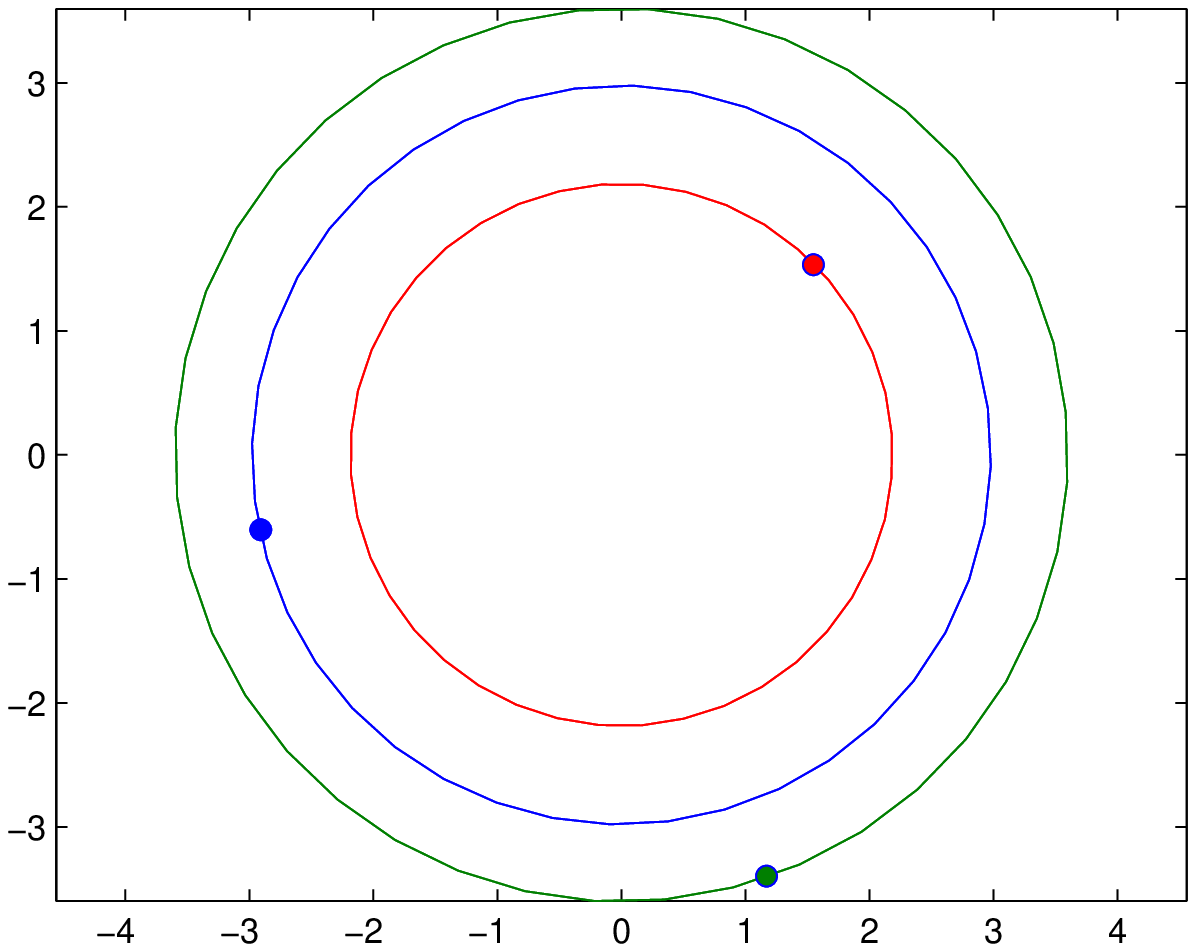}
    \parbox{2in}{\caption{$m=[1,0.5,1.5]$, $J=0.832$}\label{P3_N_1_040_0208}}
    \end{minipage}
\end{figure}

The following orbits are also found for the planar N-body problems,
where $m$ shows the number and the masses of the bodies, $J$ is the
corresponding value of Lagrangian functional action.
\begin{figure}[H]
    \begin{minipage}[t]{0.5\linewidth}
    \centering
    \includegraphics[width=2.7in]{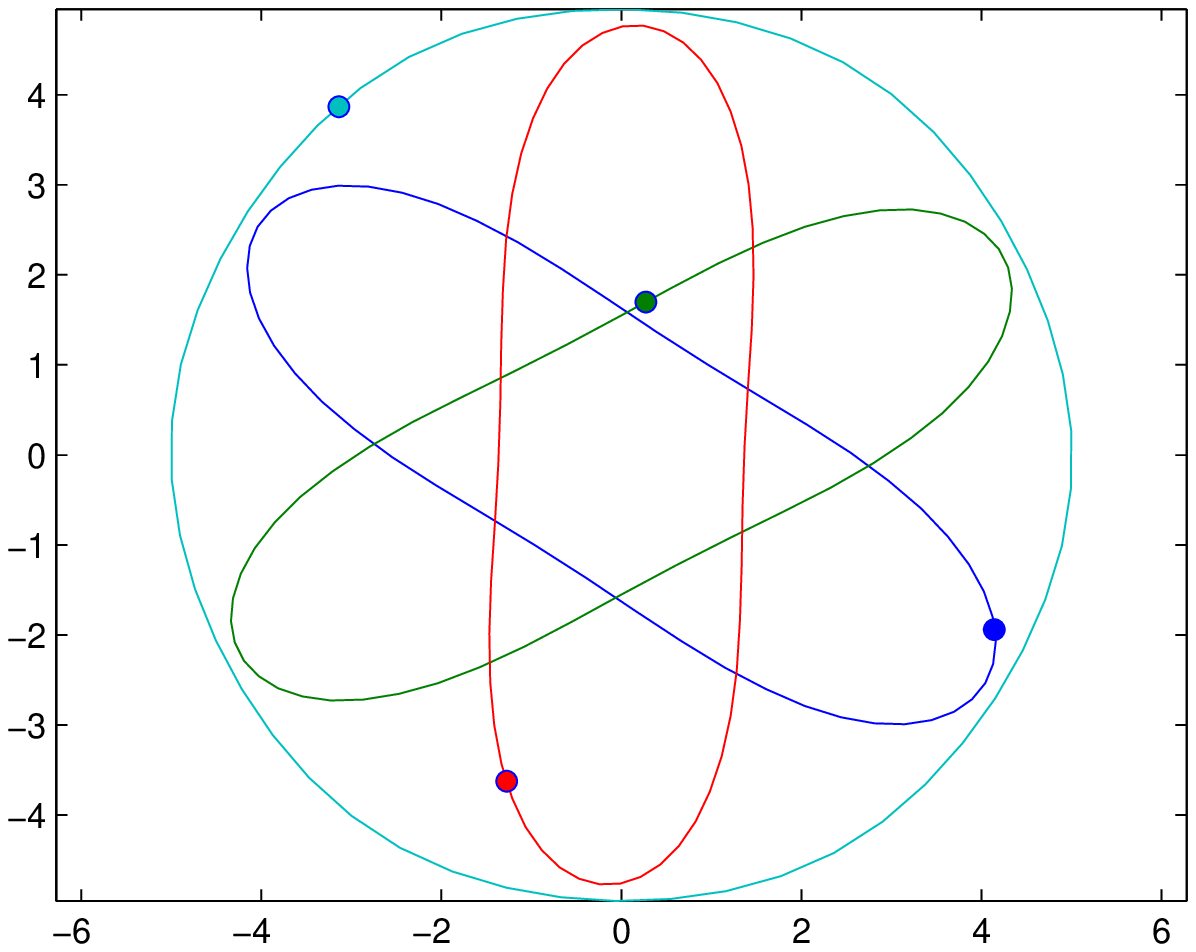}
    \parbox{2in}{\caption{$m=[1,1,1,1]$, $J=1.522$}}\label{P4_E_2_050_0127}
    \end{minipage}
    \begin{minipage}[t]{0.5\linewidth}
    \centering
    \includegraphics[width=2.7in]{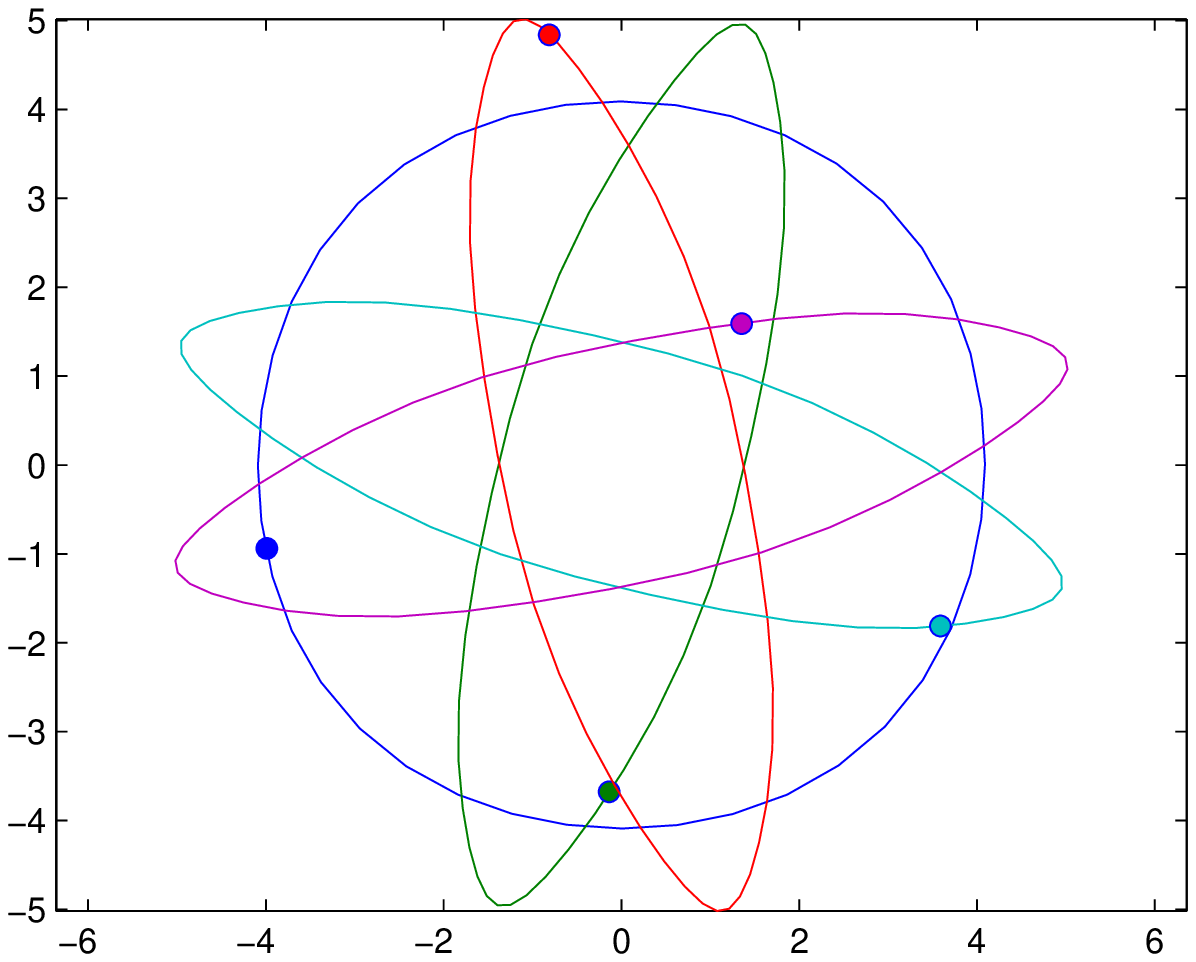}
    \parbox{2.2in}{\caption{$m=[1,1,1,1,1]$, $J=2.794$}}\label{P5_E_1_040_0127}
    \end{minipage}
\end{figure}

\begin{figure}[H]
    \begin{minipage}[t]{0.5\linewidth}
    \centering
    \includegraphics[width=2.7in]{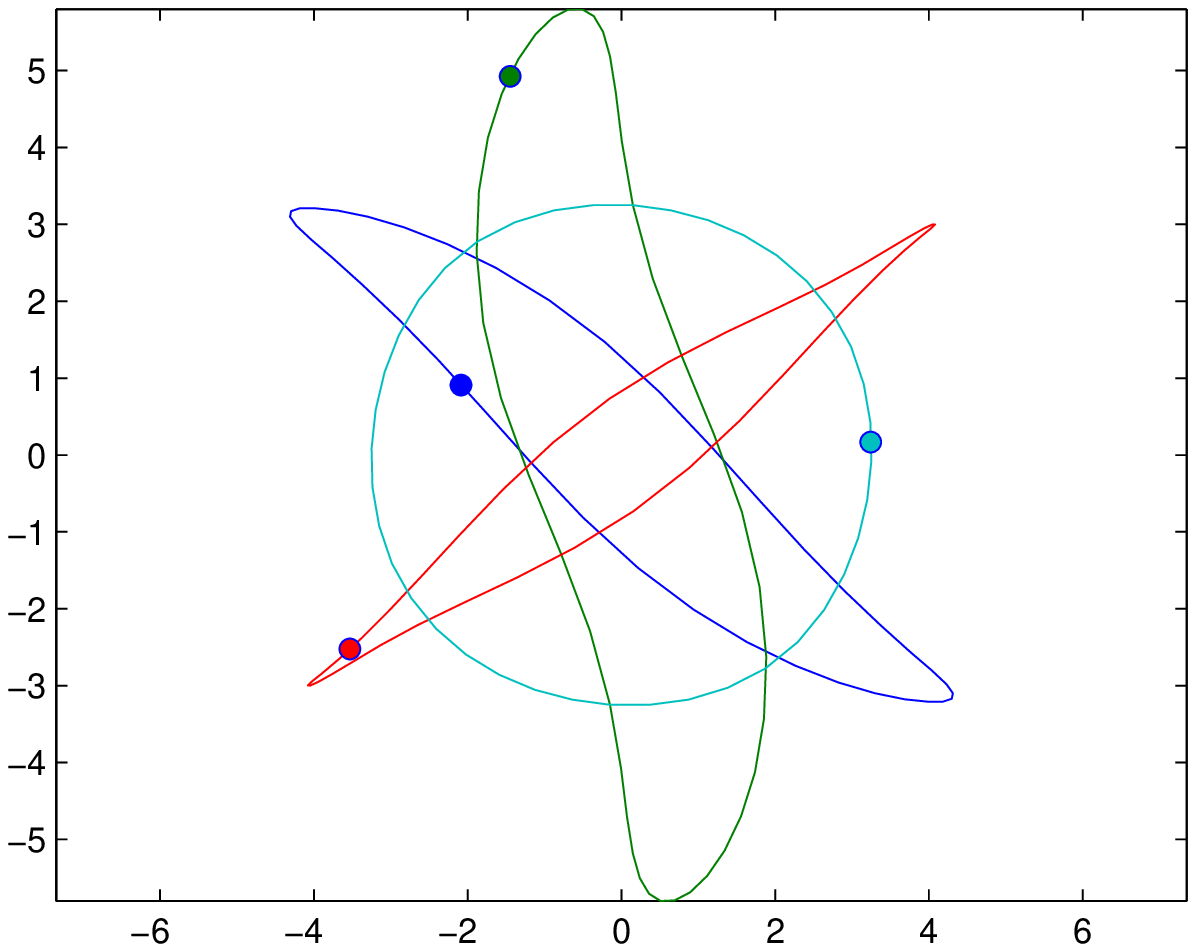}
    \parbox{2.2in}{\caption{$m=[1,0.5,1.5,2.5]$, $J=2.747$}}\label{P4_N_1_040_0209}
    \end{minipage}
    \begin{minipage}[t]{0.5\linewidth}
    \centering
    \includegraphics[width=2.7in]{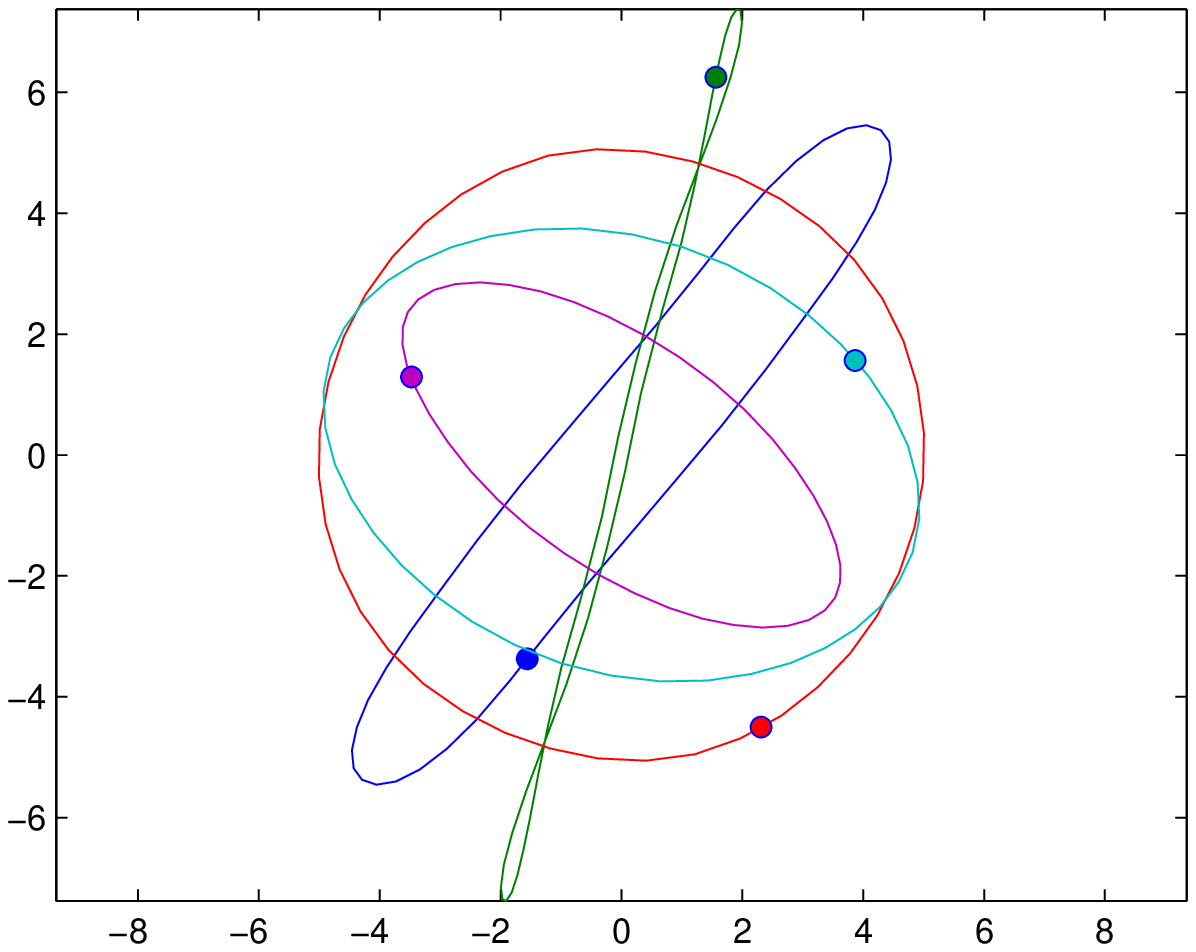}
    \parbox{2.3in}{\caption{$m=[1,0.5,1.5,2,3]$, $J=5.541$}}\label{P5_N_1_040_0209}
    \end{minipage}
\end{figure}

\begin{figure}[H]
    \begin{minipage}[t]{0.5\linewidth}
    \centering
    \includegraphics[width=2.7in]{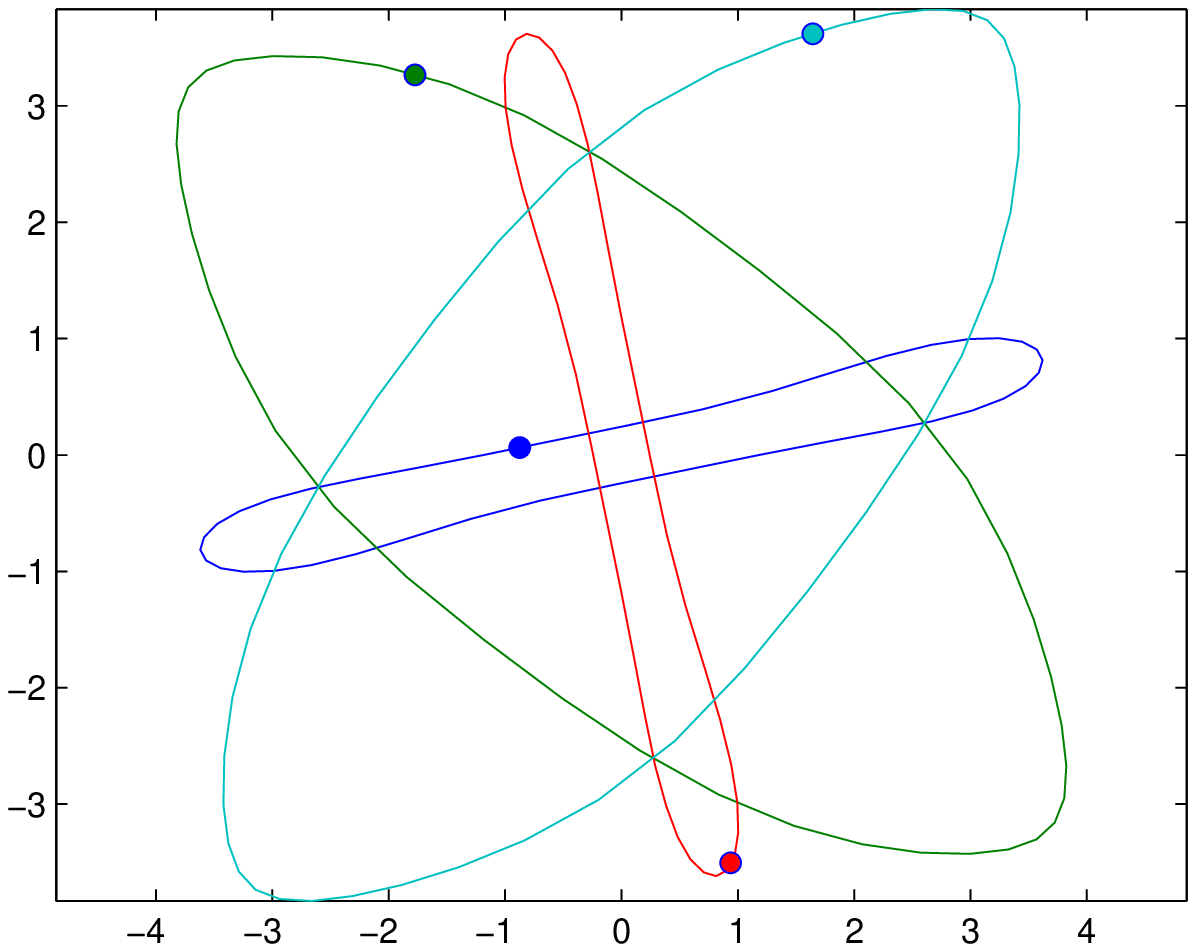}
    \parbox{2.2in}{\caption{$m=[1,0.5,1,0.5]$, $J=1.086$}}\label{P4_N_2_040_0330}
    \end{minipage}
    \begin{minipage}[t]{0.5\linewidth}
    \centering
    \includegraphics[width=2.7in]{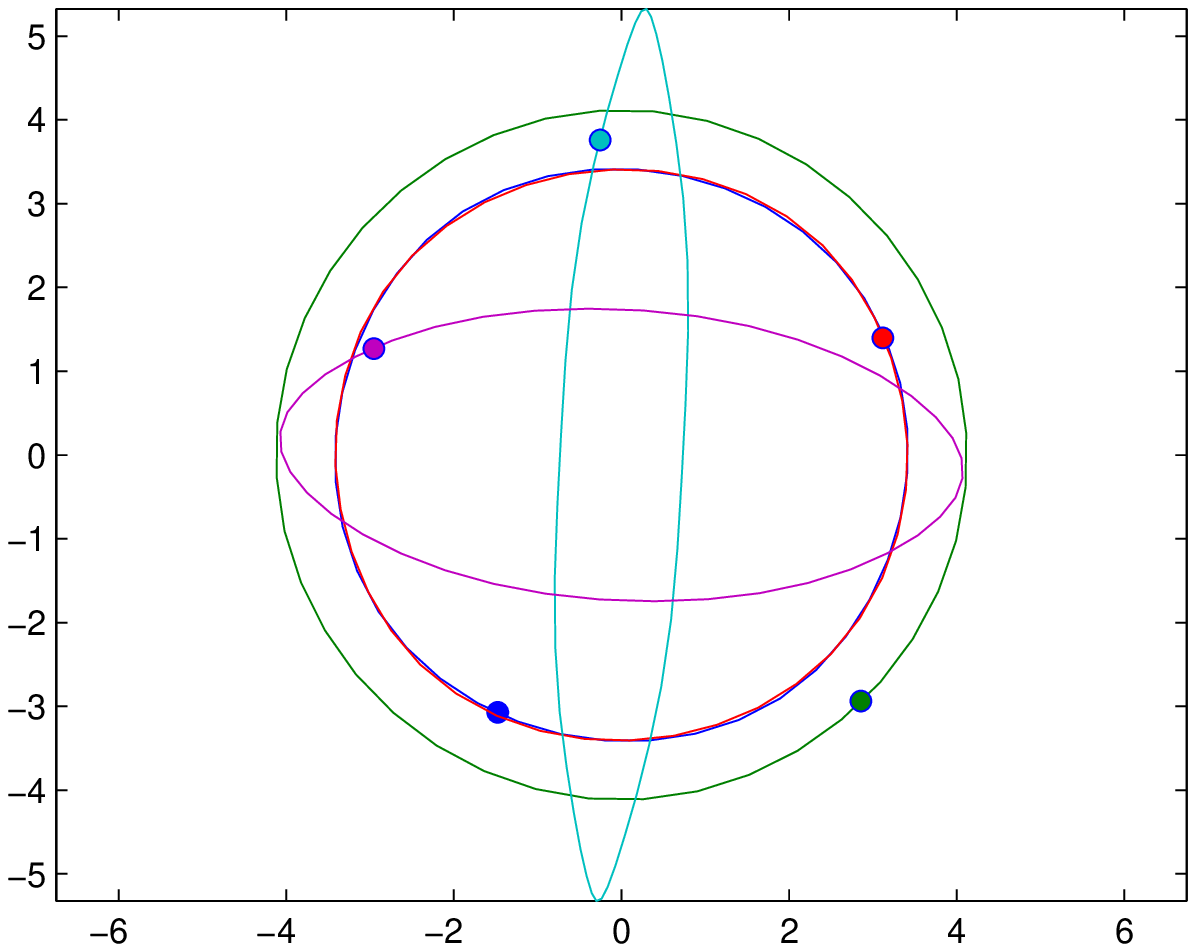}
    \parbox{2.3in}{\caption{$m=[1,0.5,1,0.5,1]$, $J=1.827$}}\label{P5_N_3_040_0508}
    \end{minipage}
\end{figure}

For spacial N-body problems,
\begin{equation*}
q^{(1,2,3)}\left(t+\frac{T}{2}\right)=-q(t)^{(1,2,3)},\quad
q^{(1,2,3)}\in R^{3},
\end{equation*}
i.e. $A_{1}=A_{2}=A_{3}=2$ in (\ref{Ms}), it is also a simple radial
symmetric constraint but in space. We find the orbits in the
following figures. In each figure, the picture on the top left
corner is the projection of the orbits on $x$-$y$ plane, the picture
on the top right corner is the projection of the orbits on $y$-$z$
plane, the picture on the bottom left corner is the projection of
the orbits on $z$-$x$ plane, and the picture on the bottom right
corner is the orbits in space.
\begin{figure}[H]
    \begin{minipage}[t]{0.5\linewidth}
    \centering
    \includegraphics[width=3in]{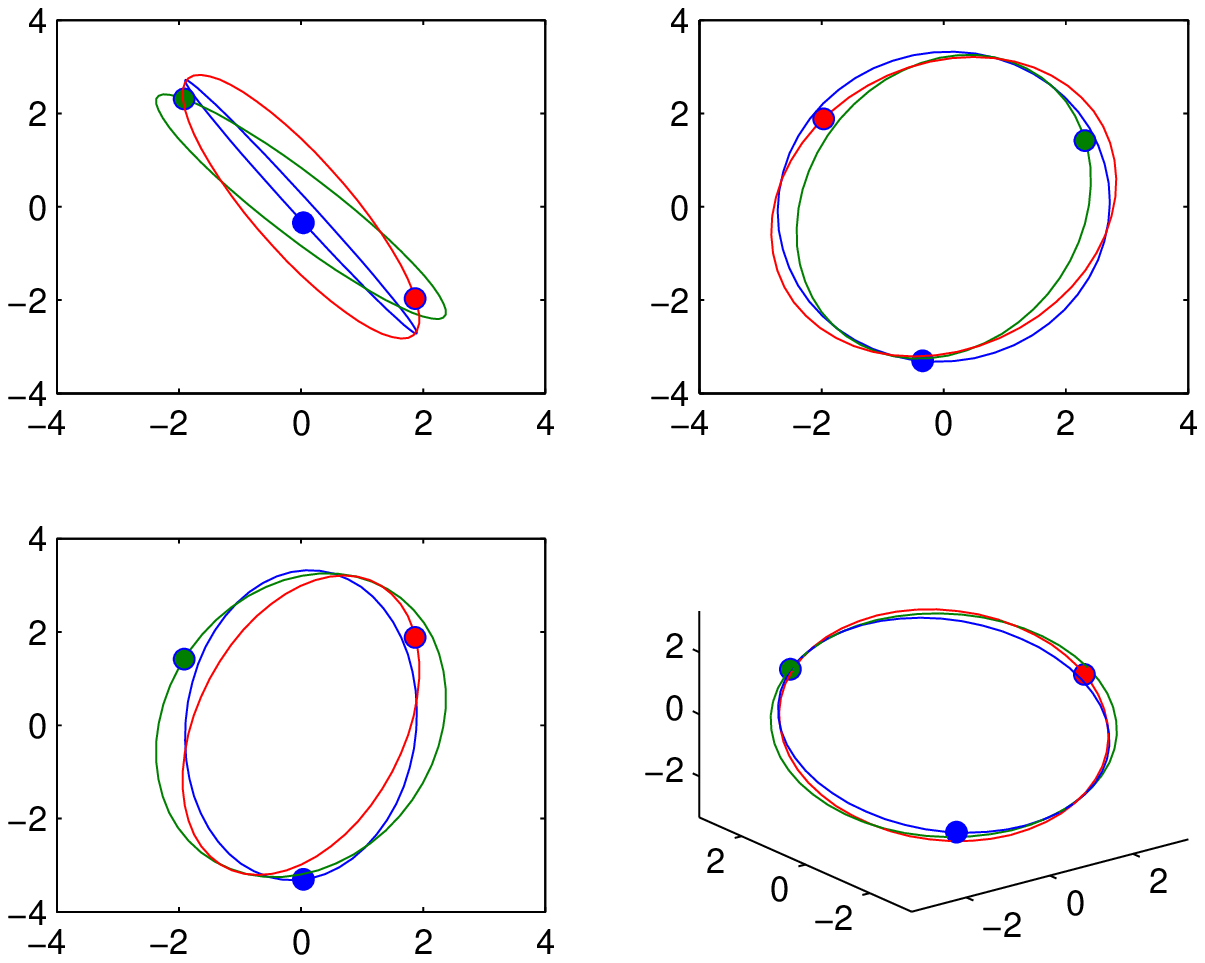}
    \parbox{2.2in}{\caption{$m=[1,1,1]$, $J=0.783$}}\label{S3_E_1_050_0126}
    \end{minipage}
    \begin{minipage}[t]{0.5\linewidth}
    \centering
    \includegraphics[width=3in]{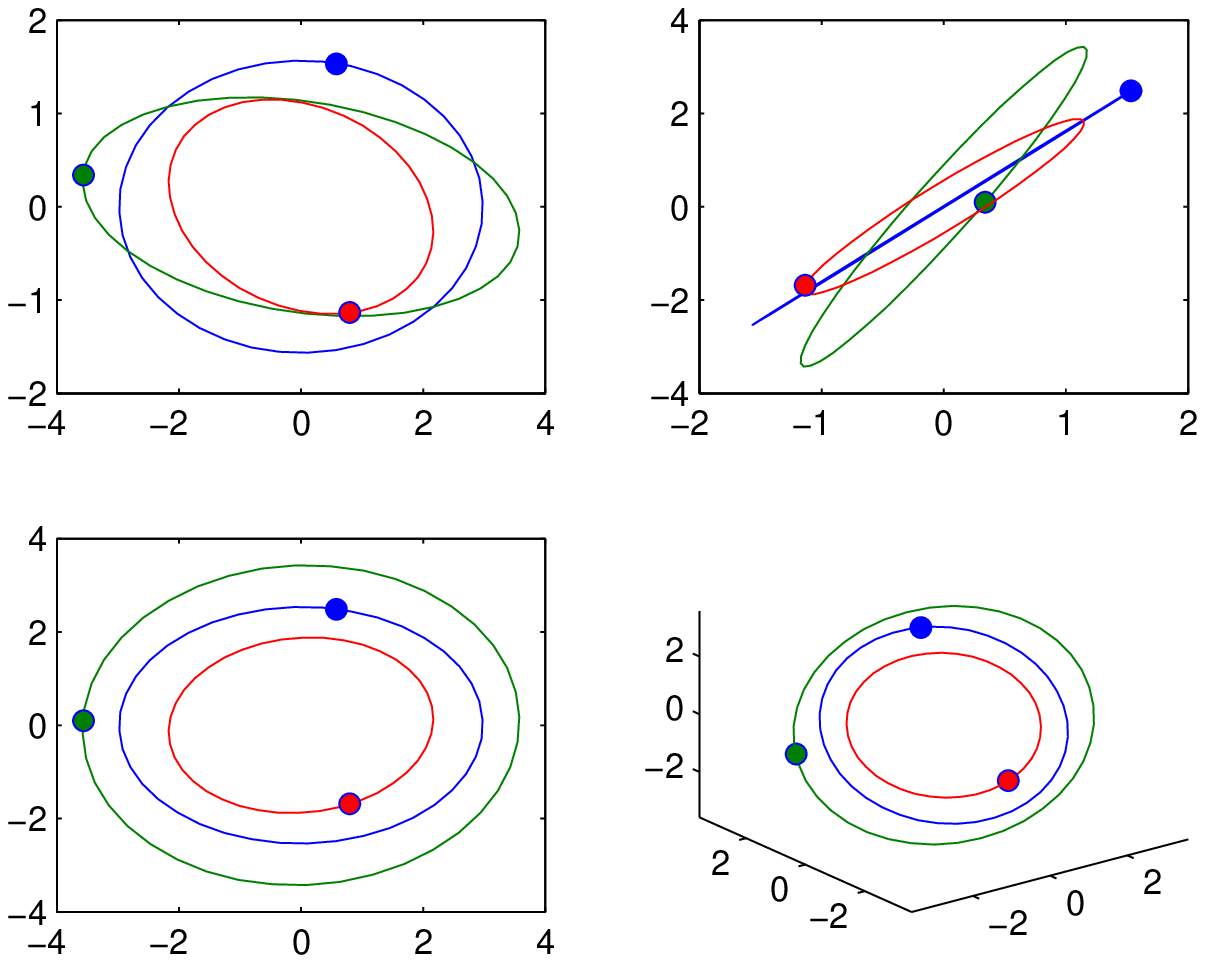}
    \parbox{2.2in}{\caption{$m=[1,0.5,1.5]$, $J=0.832$}}\label{S3_N_1_040_0209}
    \end{minipage}
\end{figure}

\begin{figure}[H]
    \begin{minipage}[t]{0.5\linewidth}
    \centering
    \includegraphics[width=3in]{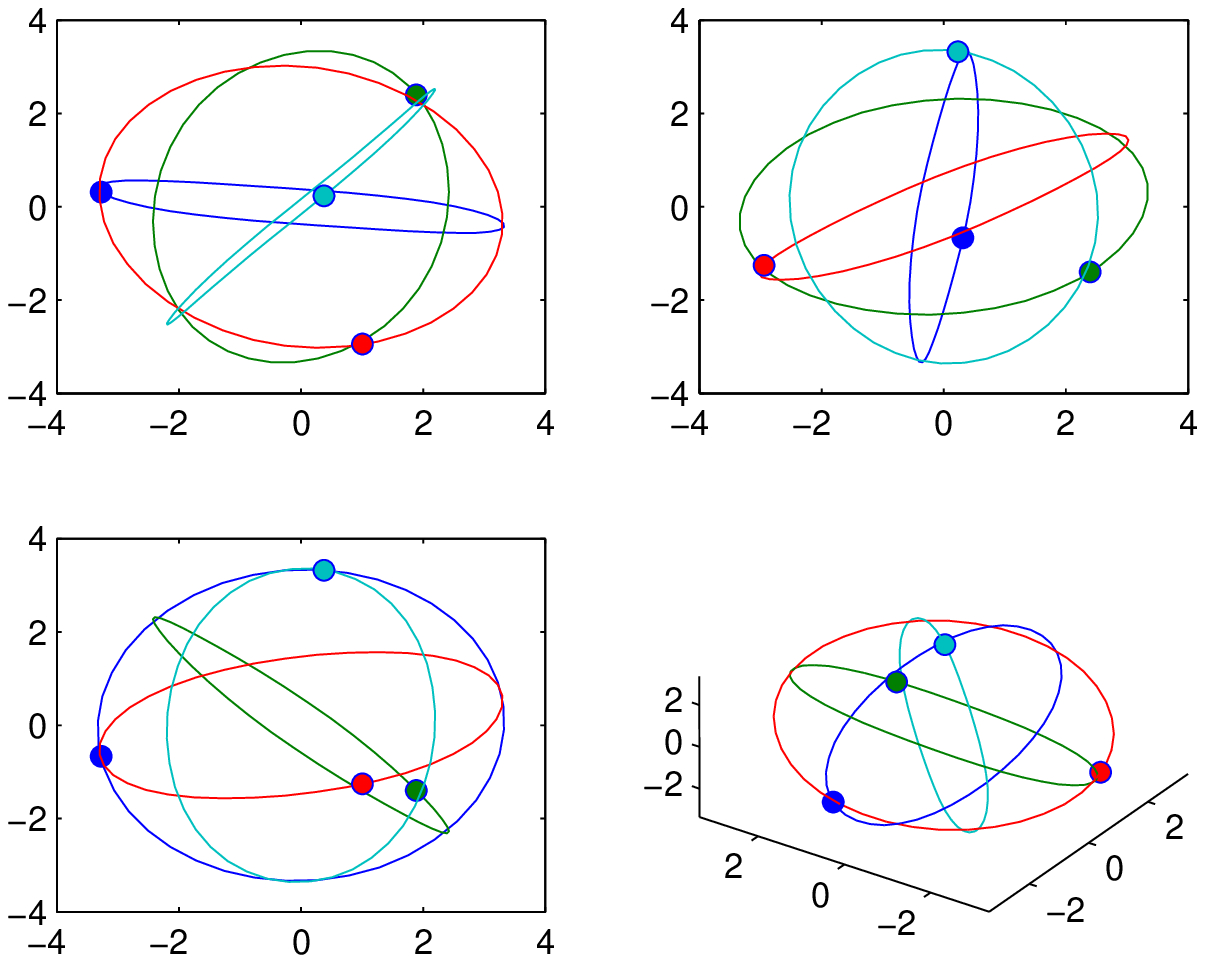}
    \parbox{2.2in}{\caption{$m=[1,1,1,1]$, $J=1.665$}}\label{S4_E_1_040_0202}
    \end{minipage}
    \begin{minipage}[t]{0.5\linewidth}
    \centering
    \includegraphics[width=3in]{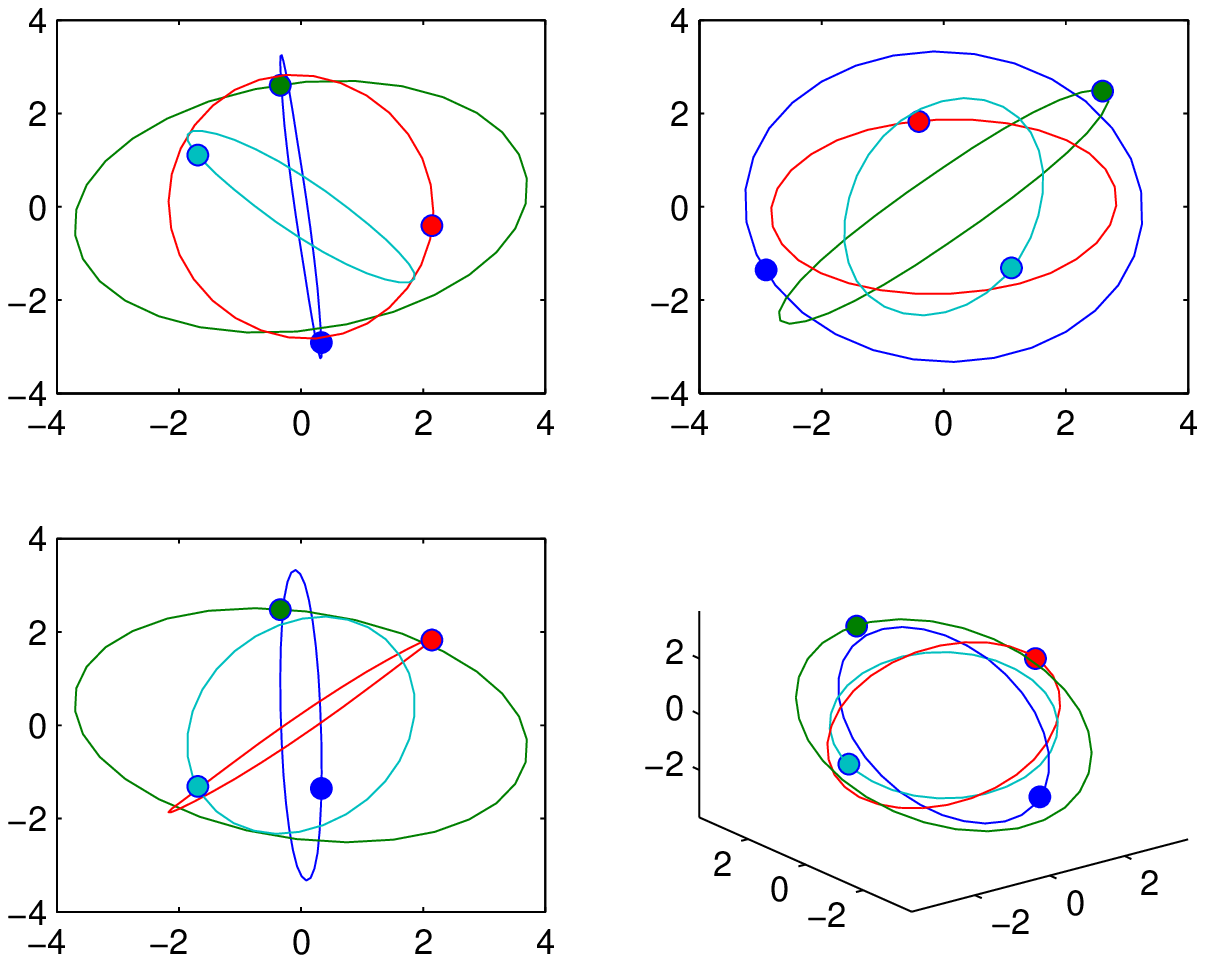}
    \parbox{2.2in}{\caption{$m=[1,0.5,1.5,2]$, $J=2.721$}}\label{S4_N_1_030_0208}
    \end{minipage}
\end{figure}

\begin{figure}[H]
    \begin{minipage}[t]{0.5\linewidth}
    \centering
    \includegraphics[width=3in]{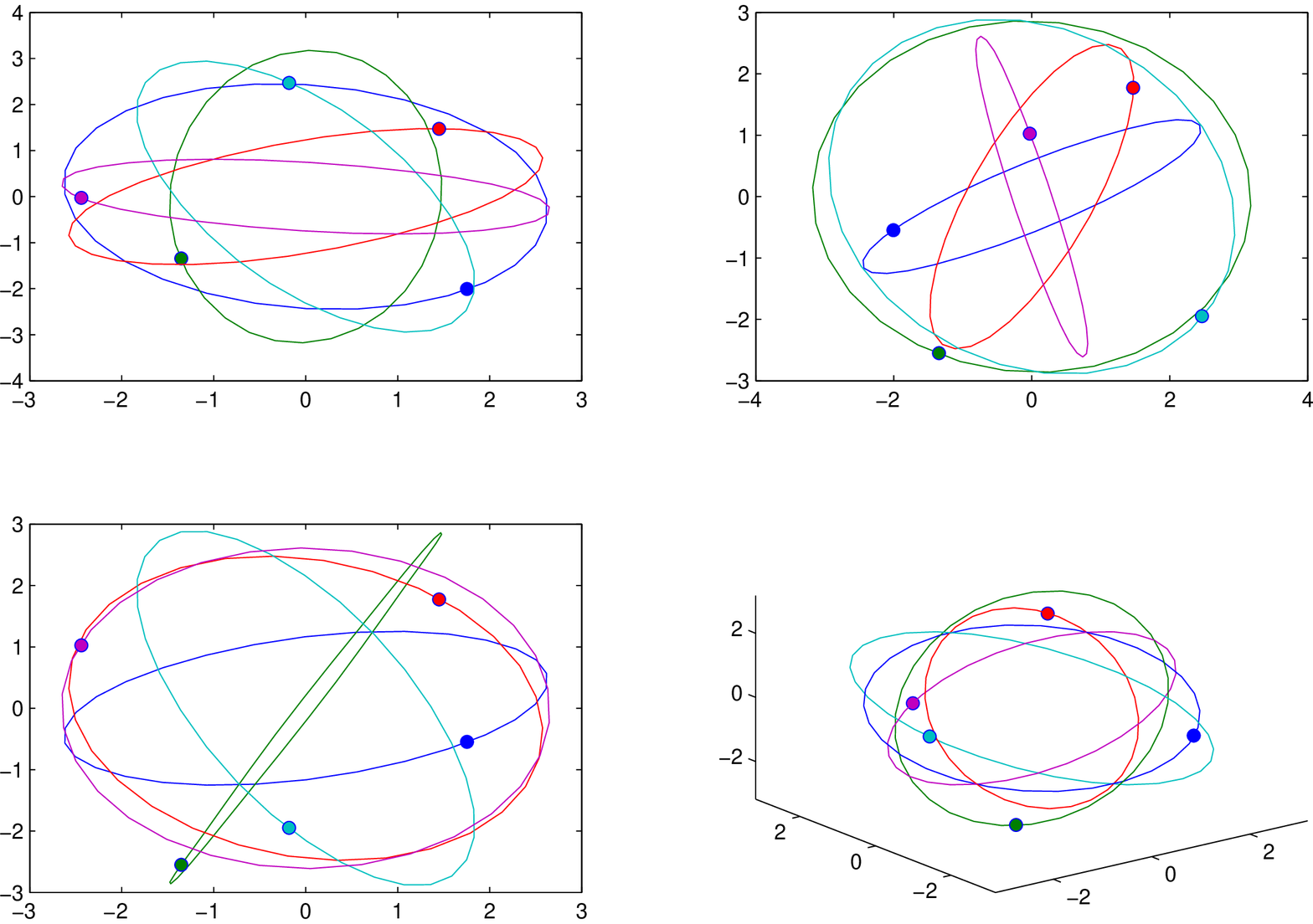}
    \parbox{2.3in}{\caption{$m=[1,0.5,1,0.5,1]$, $J=2.112$}}\label{S5_N_2_040_0330}
    \end{minipage}
    \begin{minipage}[t]{0.5\linewidth}
    \centering
    \includegraphics[width=3in]{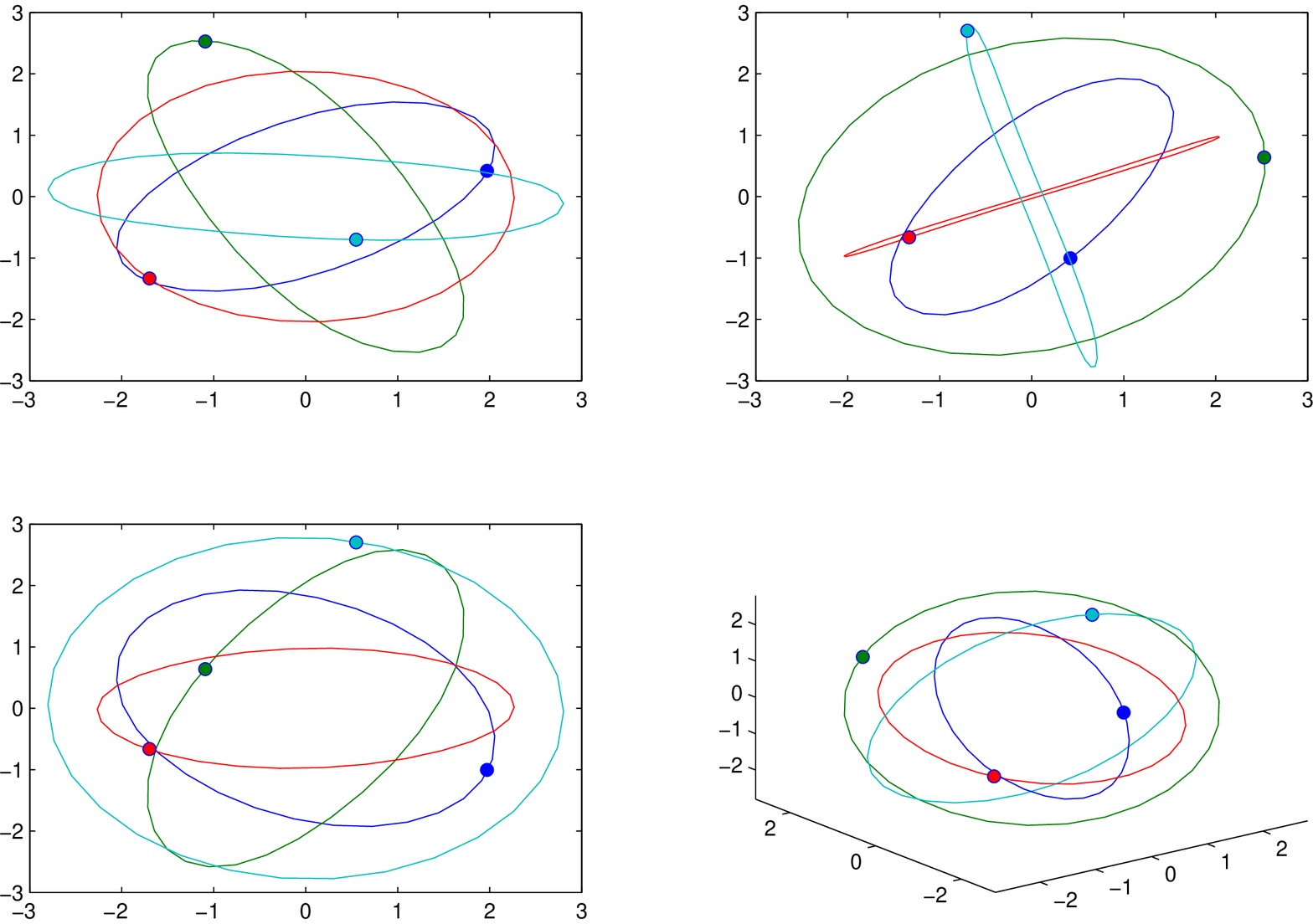}
    \parbox{2.2in}{\caption{$m=[1,0.5,1,0.5]$, $J=1.200$}}\label{S4_N_2_030_0330}
    \end{minipage}
\end{figure}

\begin{figure}[H]
    \begin{minipage}[t]{0.5\linewidth}
    \centering
    \includegraphics[width=3in]{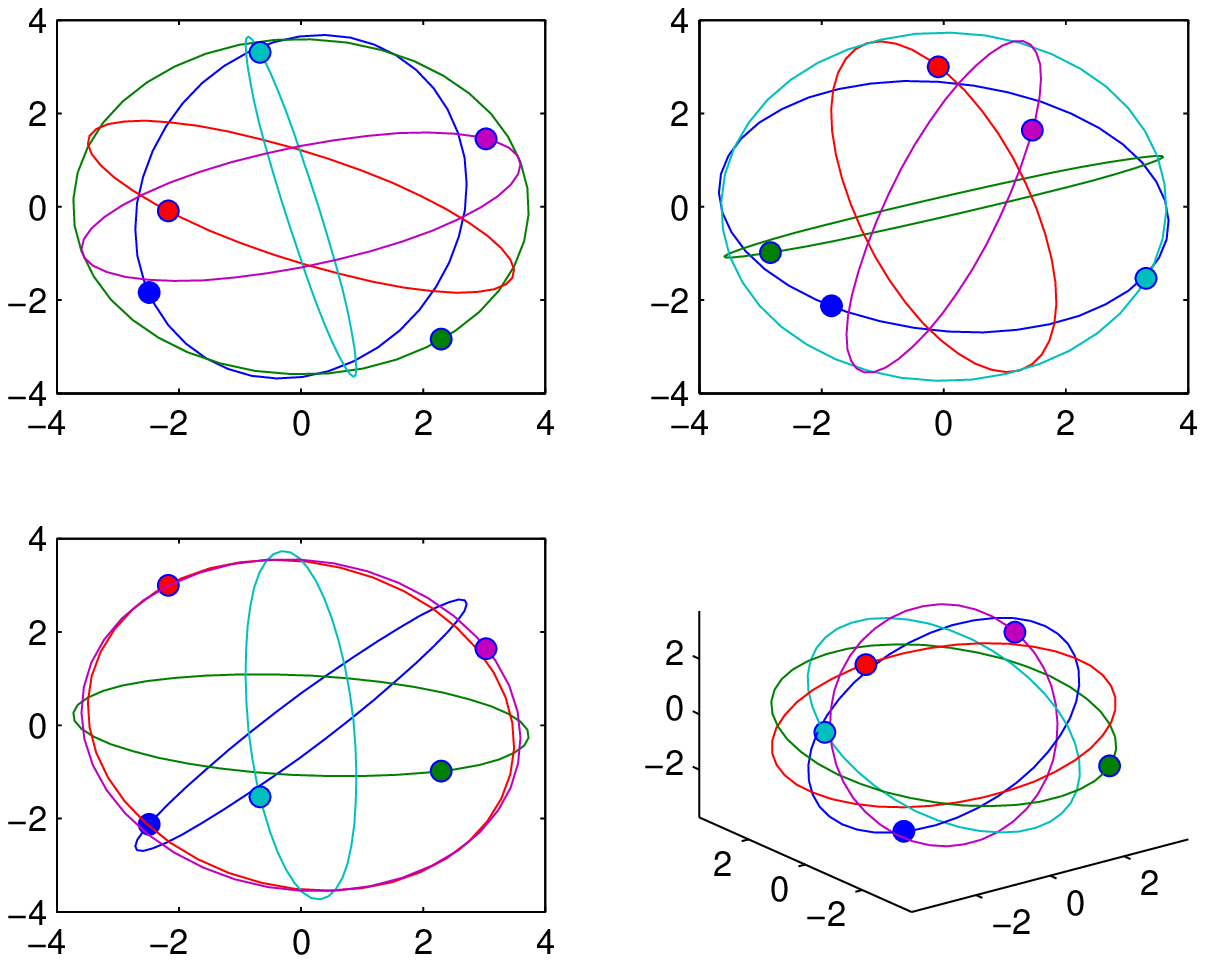}
    \parbox{2.2in}{\caption{$m=[1,1,1,1,1]$, $J=2.596$}}\label{S5_E_1_040_0209}
    \end{minipage}
    \begin{minipage}[t]{0.5\linewidth}
    \centering
    \includegraphics[width=3in]{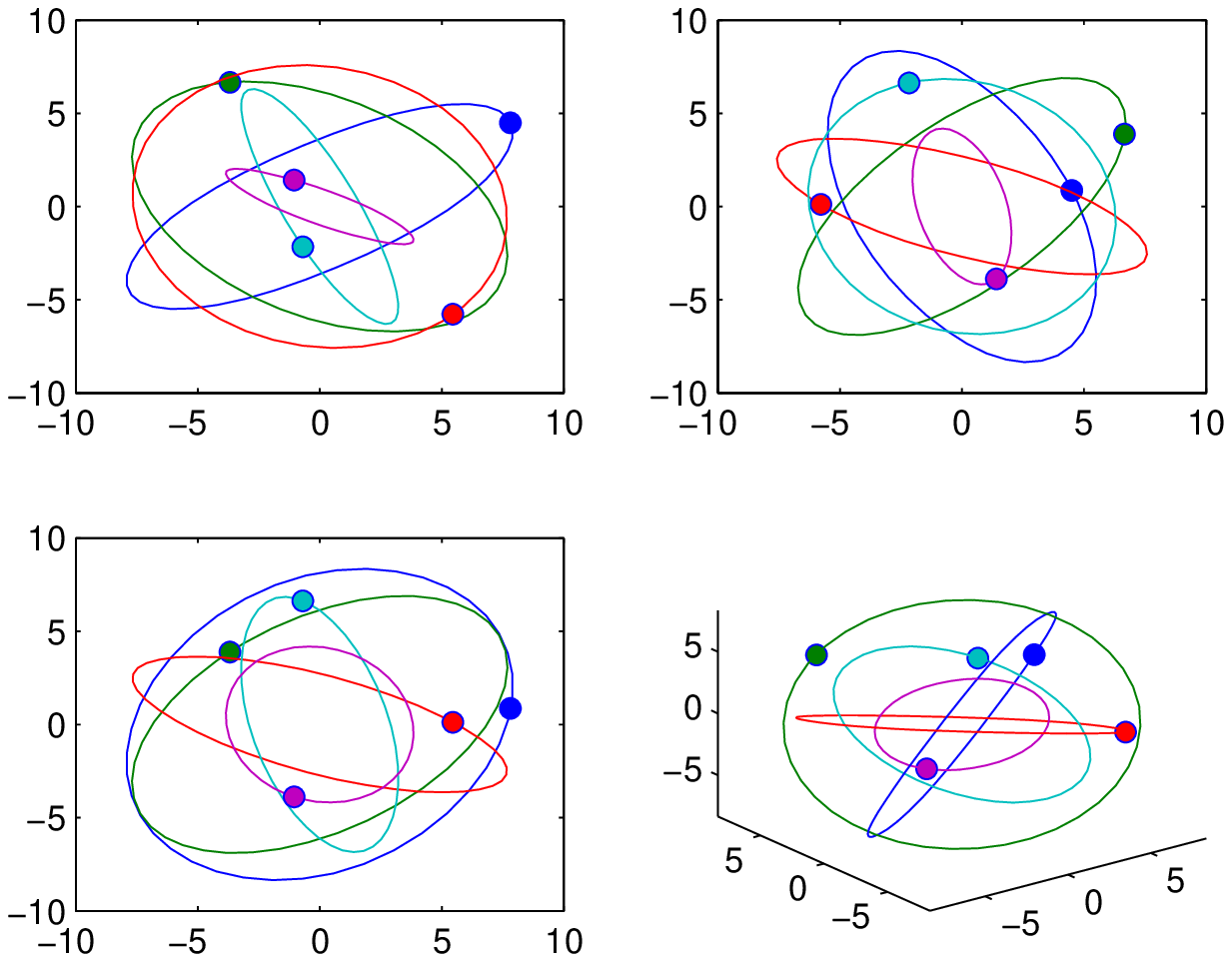}
    \parbox{2.2in}{\caption{$m=[1,2,4,8,16]$, $J=43.348$}}\label{S5_N_1_040_0209}
    \end{minipage}
\end{figure}

When we set $A_{1}=4, A_{2}=2$ in (\ref{Ms}),
$R^{d}=R^{3}=R^{2}\times R^{1}$,  we have the following multi-radial
symmetric constraint
\begin{equation*}
q_{i}^{(1,2)}(t+ \frac{T}{4}) = -q_{i}^{(1,2)}(t), \qquad
q_{i}^{(1,2)}\in R^{2},
\end{equation*}
\begin{equation*}
q_{i}^{(3)}(t+ \frac{T}{2}) = -q_{i}^{(3)}(t), \qquad q_{i}^{(3)}\in
R^{1}.
\end{equation*}
In this case, the following orbits are found.

\begin{figure}[H]
    \begin{minipage}[t]{0.5\linewidth}
    \centering
    \includegraphics[width=3in]{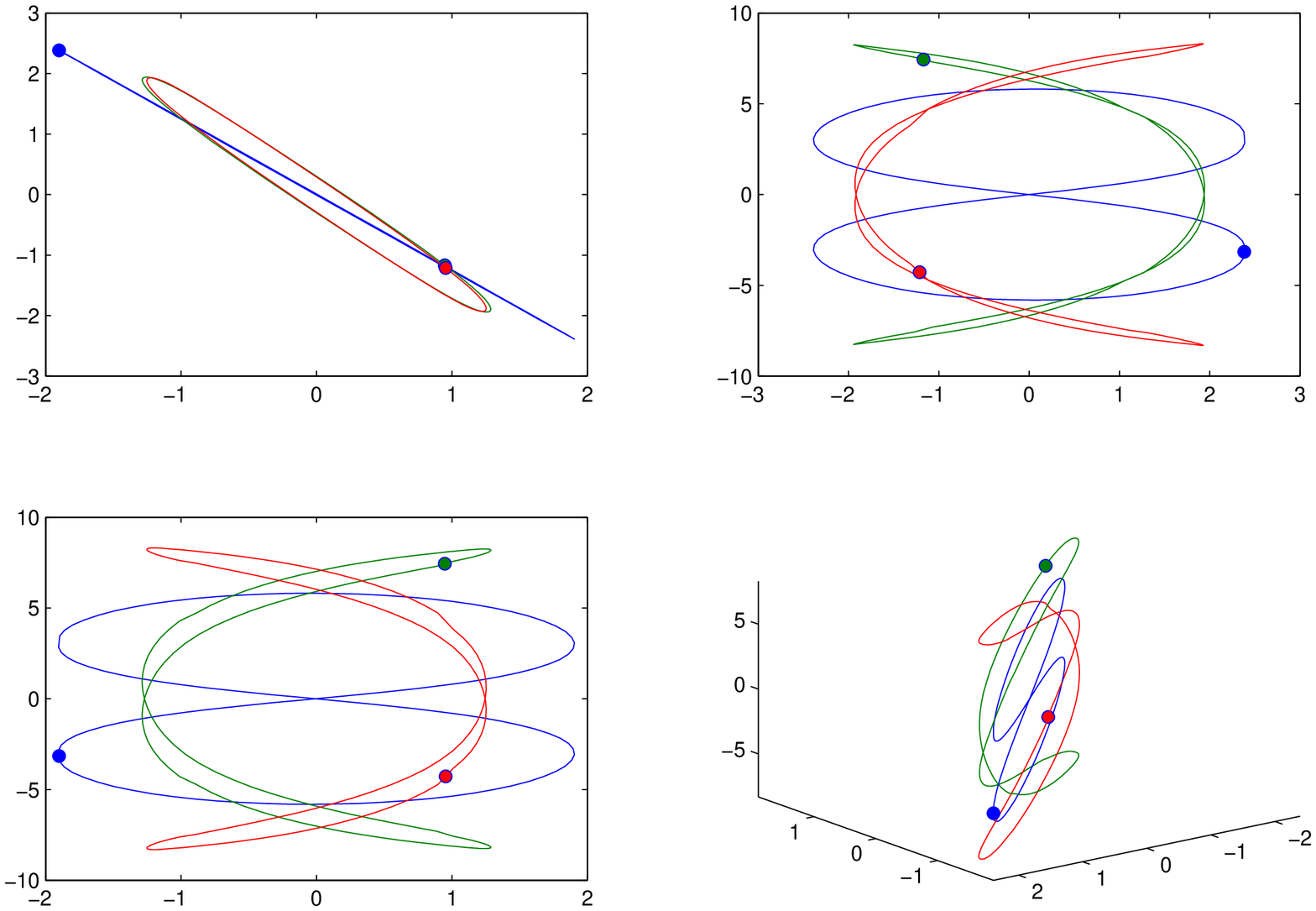}
    \parbox{2.2in}{\caption{$m=[1,1,1]$, $J=1.097$}}\label{S3_E_1_sym4T_30pt_0330}
    \end{minipage}
    \begin{minipage}[t]{0.5\linewidth}
    \centering
    \includegraphics[width=3in]{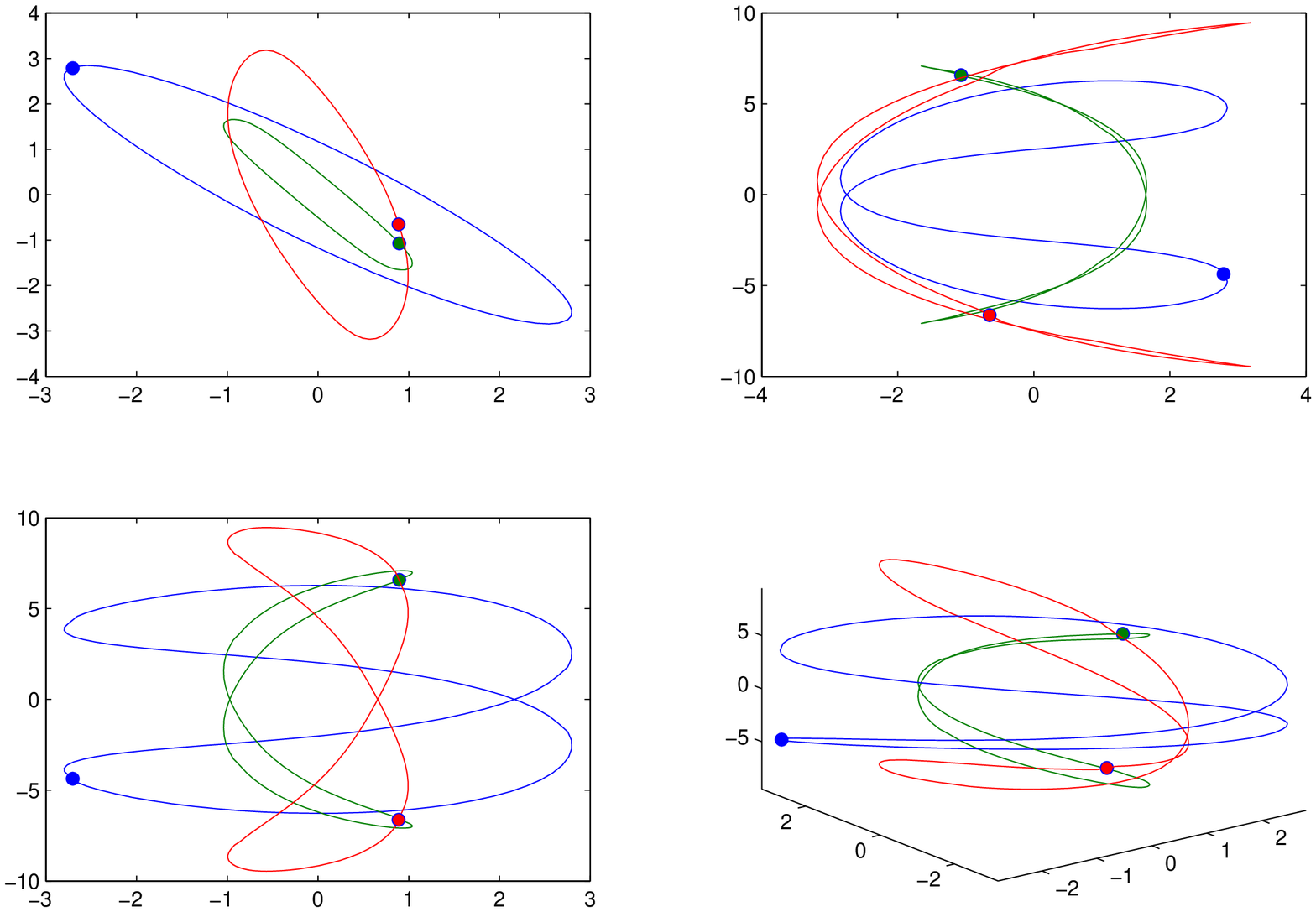}
    \parbox{2.2in}{\caption{$m=[1,2,1]$, $J=0.839$}}\label{S3_N_1_sym4T_30pt_0330}
    \end{minipage}
\end{figure}

\begin{figure}[H]
    \begin{minipage}[t]{0.5\linewidth}
    \centering
    \includegraphics[width=3in]{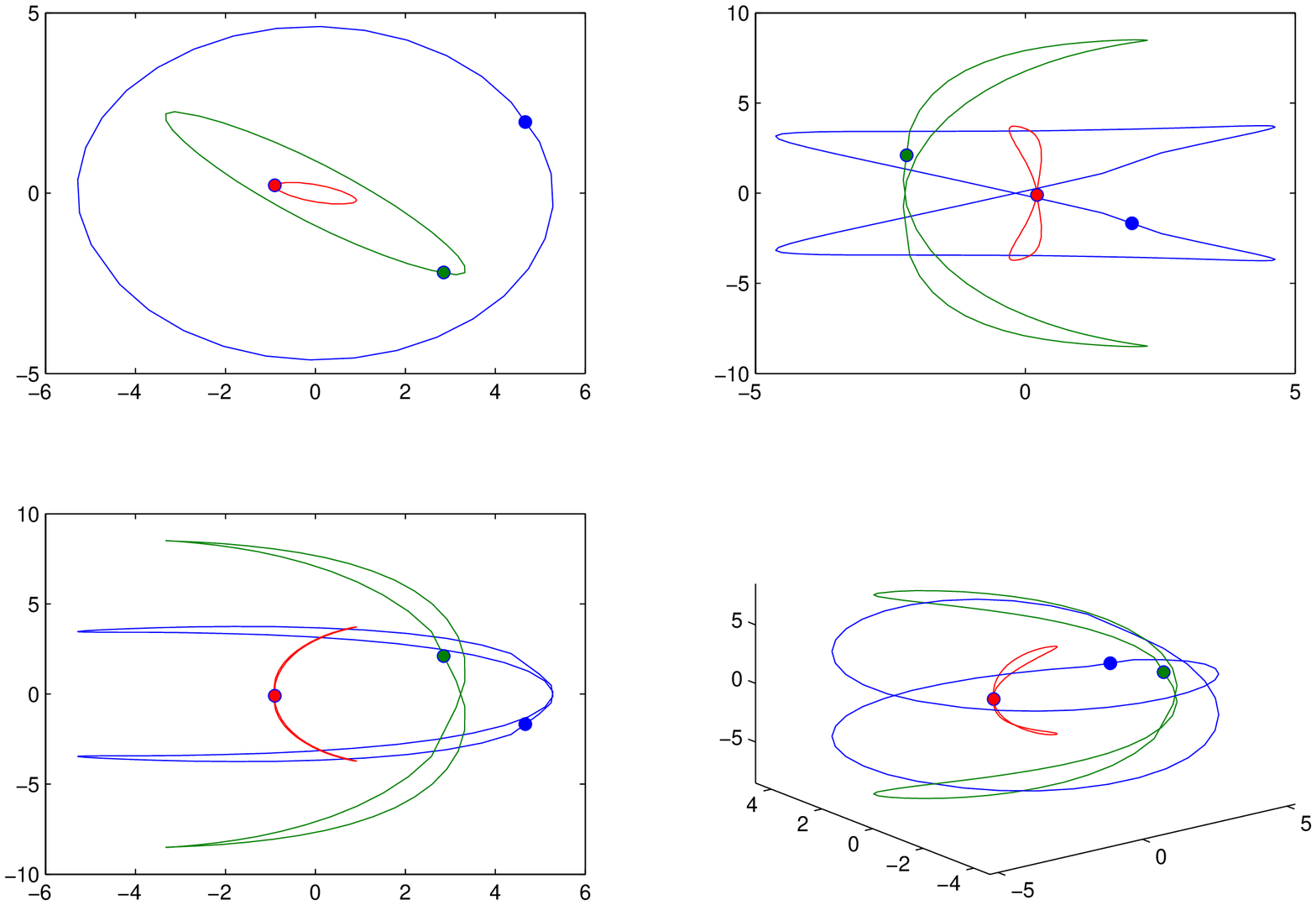}
    \parbox{2.2in}{\caption{$m=[1,2,11]$, $J=5.716$}}\label{Multi_N_3_192_060730}
    \end{minipage}
    \begin{minipage}[t]{0.5\linewidth}
    \centering
    \includegraphics[width=3in]{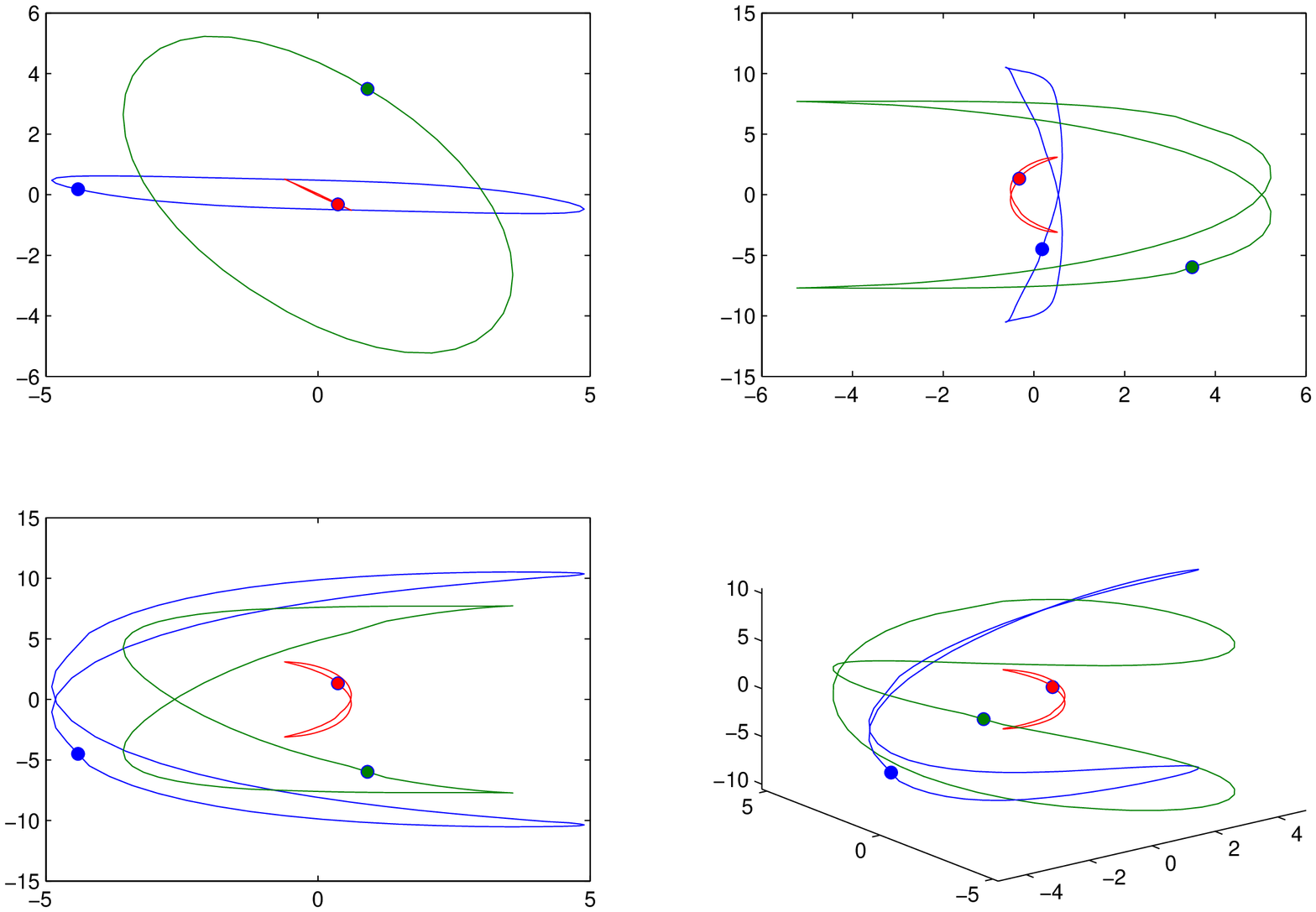}
    \parbox{2.2in}{\caption{$m=[1,1,10]$, $J=3.216$}}\label{Multi_N_3_240_060730}
    \end{minipage}
\end{figure}

\section{Conclusion}\label{sec:5}

The N-body type difference equation is first introduced in this
paper with a new difference scheme. To prove the existence of
periodic solutions for N-body type problem, the variational approach
and the method of minimizing the Lagrangian functional action, both
of which are typical in studying N-body dynamical systems, are
employed in this paper. In finding the generalized solutions and
non-collision solutions numerically, we use the method of minimizing
the functional action of N-body type difference equation. The
presented numerical solutions of N-body problem illuminate that the
variational difference method, which is combined with the
traditional steepest method, is rather effective in finding the
periodic solutions for N-body type problem numerically. with these
numerical examples, we can conclude that there exist some
interesting orbits under multi-radial symmetric constraint; and
surprisingly, even under radial symmetric constraint we have
obtained some interesting orbits. Furthermore, difference equations
transform infinity dimension dynamical systems into finite dimension
dynamical systems, and it can be used in the studies for other
solutions under corresponding symmetric or choreographic
constraints.

\end{document}